\documentclass[a4paper,11pt]{article}

\usepackage[affil-it]{authblk}
\usepackage{amsmath}
\usepackage{amsfonts}
\usepackage{amssymb}
\usepackage{graphicx}
\DeclareGraphicsExtensions{.jpg,.png,.pdf,.eps,.bmp,.gif}
\usepackage[latin1]{inputenc}
\title{A white noise approach to stochastic integration with respect to the Rosenblatt process.}
\author{Benjamin Arras
\thanks{Electronic address: \texttt{benjamin.arras@ecp.fr}}}
\affil{Ecole Centrale Paris and INRIA Regularity team\\ 
Grande Voie des Vignes, 92295~Chatenay-Malabry, France}

\begin{document}
\maketitle

\begin{abstract}
In this paper, we define a stochastic calculus with respect to the Rosenblatt process by means of white noise distribution theory. For this purpose, we compute the translated characteristic function of the Rosenblatt process at time $t>0$ in any direction $\xi\in S(\mathbb{R})$ and the derivative of the Rosenblatt process in the white noise sense. Using Wick multiplication by the former derivative and Pettis integration, we define our stochastic integral with respect to the Rosenblatt process for a wide class of distribution processes. We obtain an explicit formula for the variance of such a stochastic integral and Itô formulae for a certain class of functionals of the Rosenblatt process. Finally, we compare our stochastic integral to other approaches.
\end{abstract}

\\
Since the nineties, there have been mainly two different approaches in order to develop stochastic calculus for processes which are not semimartingales. The first one is a pathwise Riemann-Stieltjes integration which originates from a stochastic analysis of fractional Brownian motion for $H>\frac{1}{2}$ (see \cite{lin1995stochastic}). Indeed, fractional Brownian motion (fBm), popularized by Mandelbrot (see \cite{MR0242239}), is not a semimartingale for $H\ne\frac{1}{2}$. This approach which has been extended by Zhäle, using fractional integration by part formula, (see \cite{zahle1998integration}), allows us to obtain an integral with similar properties to the ones of the Stratonovich integral. The second approach is based on tools from the Malliavin calculus and leads us to a stochastic integral close to the Itô one. These techniques were first developped by Decreusefond and Üstünel (see \cite{MR1677455}) for fractional Brownian motion and extended by Alos et al. to Gaussian Volterra processes (see \cite{alos2001stochastic}). In \cite{MR1956473}, using white noise distribution theory, the author developped a stochastic calculus with respect to fractional Brownian motion (for every $H\in(0,1)$) which leads, quite easily, to an Itô formula for generalized functionals of this process.\\
\\
Moreover, fBm belongs to the family of Hermite processes which appear in non-central limit theorems for processes defined as integrals or partial sums of non-linear functionals of stationary Gaussian sequences with long-range dependence (see \cite{MR550123}, \cite{taqqu1975weak} and \cite{MR550122}). They admit the following representation, for all $d\geq 1$:
\begin{align*}
\forall t>0\quad Y^{H,d}_t=c(H_0)\int_{\mathbb{R}}...\int_{\mathbb{R}}\left(\int_0^t\prod_{j=1}^d(s-x_j)^{H_0-1}_+ds\right)dB_{x_1}...dB_{x_d}
\end{align*}
where $\{B_x:\ x\in\mathbb{R}\}$ is a two-sided Brownian motion, $c(H_0)$ is a normalizing constant such that $\mathbb{E}[|Y^{H,d}_1|^2]=1$ and $H_0=\frac{1}{2}+\frac{H-1}{d}$ with $H\in(\frac{1}{2},1)$. For $d=1$, one recovers fractional Brownian motion. These processes share many properties with fBm. Namely, they are $H$-self-similar processes with stationary increments. They possess the same covariance structure, exhibit long range-dependence and their sample paths are almost-surely $\delta$-Hölder continuous, for every $\delta<H$. For $d=2$, the process is called the Rosenblatt process by Taqqu (see \cite{taqqu1975weak}).\\
\\
This process has received lots of interest in the past and more recent years. The distributional properties have first been studied by Albin (see \cite{MR1650532}) and more recently by Veillette and Taqqu (\cite{Veillette}) and Maejima and Tudor (\cite{Maejima20131490}). A wavelet-type expansion has been established by Pipiras (\cite{MR2105535}). Geometric properties of random fractal sets generated by anisotropic multidimensional version of the Rosenblatt process have been considered in \cite{MR2759163}. The construction of estimators for the self-similarity parameter $H$ has been investigated by Bardet and Tudor in \cite{Bard} and by Tudor and Viens in \cite{tudor2009variations}.\\
\\
Finally, stochastic calculus with respect to the Rosenblatt process has been developped by Tudor in \cite{MR2374640} from both, the pathwise type calculus and Malliavin calculus points of view. Even if these two approaches are successful in order to define a stochastic integral with respect to the Rosenblatt process, the Malliavin calculus one fails to give an Itô formula for the Rosenblatt process in the divergence sense for general smooth functionals. Indeed, the appearance of two trace terms from a limiting procedure makes the Itô formula in the divergence sense difficult to obtain (see Section $8$ of \cite{MR2374640} and more specifically Theorem $3$). However, by analyzing the peculiar cases of the monomials $x^2$ and $x^3$, the author noticed the appearance of certain non-zero cumulants of the Rosenblatt distribution (see section $2$ below for a definition) as well as the appearance of the third derivative of $x^3$. Moreover, in Remark $8$ of \cite{MR2374640}, the author expected that, for sufficiently smooth functionals, Itô formula regarding stochastic integration with respect to the Rosenblatt process should involve all the derivatives of such functionals as well as all the non-zero cumulants of the Rosenblatt distribution. Using white noise distribution theory, we define, in this article, a stochastic integral with respect to the Rosenblatt process which leads to an Itô formula for a certain class of functionals of this process confirming the above conjecture. This class of functionals is completely characterized by Paley-Wiener-Schwartz theorems. All the derivatives of these functionals appear in the change of variable formula (see section 3, specifically Remarks $3.17$ and $3.19$). The non-zero cumulants of the Rosenblatt distribution seem to be the ones responsible for this result (see Theorems $3.18$ and $3.20$).\\
\\
Last but not least, we obtain in section $3$ of this article (see Theorem $3.13$) an explicit formula for the variance of the stochastic integral with respect to the Rosenblatt process. We note the appearance of three terms. The first one is directly linked with the covariance structure of the Rosenblatt process. Since fractional Brownian motion and the Rosenblatt process share the same covariance function, this term is similar to the one appearing in the definition of the class of deterministic integrands for which one can define a Wiener integral with respect to fBm (see \cite{MR1790083}). The second one needs the use of the derivative operator from white noise distribution theory (see Definition-Theorem $1.8$ below). This term is somehow classical in this setting (see Theorem $13.16$ of \cite{MR1387829}). Finally, the third term emphasizes the non-Gaussian nature of the Rosenblatt process. Since the Rosenblatt process admits a double Wiener-Itô integral representation, the second order derivative operator is involved in the definition of this term. Moreover, this formula holds under two sufficient conditions on the Wiener-Itô expansion of the stochastic integrand process. More specifically, multidimensional fractional integrals of the kernels of the Wiener-Itô expansion are assumed to belong to a certain function space with dominating mixed smoothness (see the preamble before Theorem $3.13$ for a definition and hypothesis \textbf{H1}). These spaces were introduced by Lizorkin and Nikol'skii in \cite{MR0192223} (see \cite{MR2657119} for a recent survey of these spaces). This assumption is needed in order to define properly the second order derivative operator applied to the stochastic integrand process. The second condition ensures that one recovers the variance formula for the stochastic integral of the integrand process from the sum of the variance formulae of the stochastic integral of its projections on Wiener chaoses. This result stresses one of the advantages of white noise distribution theory over pathwise methods for the stochastic calculus with respect to the Rosenblatt process. Indeed, as in the fBm case for $H\in (1/2,1)$, one can use Young integration in order to define a stochastic integral with respect to the Rosenblatt process. Nevertheless, as stated in \cite{MR2408998} (see Remark $9$ page $17$ as well as Conclusion $3.3$ page $19$), statistical properties such as mean or variance of these stochastic integrals seem rather difficult to obtain.\\
\\
The remaining of this article is organized as follows. In the first section, we introduce the tools from the white noise distribution theory needed in order to define a stochastic calculus with respect to the Rosenblatt process. Section 2 is devoted to the study of the translated characteristic function of the Rosenblatt process at time $t>0$ in the direction $\xi\in S(\mathbb{R})$. The main results of this article are in section 3 where we define the stochastic integral by means of the Rosenblatt noise, Wick product and Pettis integration. Specifically, we derive a sufficient condition for the Rosenblatt noise integral of a stochastic process $\{\phi_t\}$ to be an element of $(L^2)$ and we compute the variance of this random variable. We obtain as well Itô formulae for functionals of the Rosenblatt process. In the last section, we compare our approach to the one in \cite{MR2374640}. Using the finite interval representation of the Rosenblatt process, we show that the Rosenblatt noise integral extends the Skorohod integral defined by Malliavin calculus tools.

\section{White noise setting.}

In this section, we briefly remind the white noise analysis introduced by Hida and al. in \cite{hida1993white}. For a good introduction to the theory of white noise, we refer the reader to the book of Kuo \cite{MR1387829}. The underlying probability space $(\Omega, \mathcal{F}, \mathbb{P})$ is the space of tempered distributions endowed with the $\sigma$-field generated by the open sets with respect to the weak* topology in $S'(\mathbb{R})$ and with the infinite dimensional Gaussian measure $\mu$ whose existence is ensured by the Bochner-Minlos theorem.\\
\\
For all $(\phi_1,...,\phi_n)\in S(\mathbb{R})$, the space of $C^{\infty}(\mathbb{R})$ functions with rapid decrease at infinity, the vector $(<.;\phi_1>,...,<.;\phi_n>)$ is a centered Gaussian random vectors with covariance matrix $(<\phi_i;\phi_j>)_{(i,j)}$. As it is written in Kuo \cite{MR1387829}, for any function $f\in L^2(\mathbb{R})$, we can define $<.;f>$ as the random variable in $L^2(\Omega, \mathcal{F}, \mathbb{P})$ obtained by a classical approximation argument and the following isometry:
$$\forall (\phi,\psi)\in S(\mathbb{R})^2\quad\mathbb{E}[<;\psi><;\phi>]=<\psi;\phi>_{L^2(\mathbb{R})}$$
Thus, for any $t\in\mathbb{R}$, we define ($\mu$-almost everywhere):
\[B_t(.)=
\begin{cases}
<.;1_{[0;t]}> & t\geq 0 \\
-<.;1_{[t;0]}> & t<0
\end{cases}
\]
From the isometry property, it follows immediately that $B_t$ is a Brownian motion on the white noise space and, by the Kolmogorov-Centsov theorem, it admits a continuous modification. Moreover, using the approximation of any function $f\in L^2(\mathbb{R})$ by step functions, we obtain:
$$<;f>=\int_{\mathbb{R}}f(s)dB_s$$
We note $\mathcal{G}$, the sigma field generated by Brownian motion and $(L^2)=L^2(\Omega,\mathcal{G},\mathbb{P})$. By the Wiener-Itô theorem, any functionals $\Phi\in (L^2)$ can be expanded uniquely into a series of multiple Wiener-Itô integrals:
$$\Phi=\sum_{n=0}^{\infty}I_n(\phi_n)$$
where $\phi_n\in \hat{L}^2(\mathbb{R}^n)$, the space of square-summable symmetric functions. Using this theorem and the second quantization operator of the harmonic oscillator operator, $A=\frac{-d^2}{dx^2}+x^2+1$, Hida and al. introduced the stochastic space of test functions $(S)$ and its dual, the space of generalized functions $(S)^*$ or Hida distributions. We refer the reader to pages 18-20 of Kuo's book for an explicit construction. We have the following Gel'fand triple:
$$(S)\subset (L^2)\subset (S)^*$$
We denote by $\left\langle \left\langle ;\right\rangle\right\rangle$ the duality bracket between elements of $(S)$ and $(S)^*$ which reduces to the classical inner product on $(L^2)$ for two elements in $(L^2)$.\\
In the context of white noise analysis, the main tool is the $S$-transform. It is a functional on $S(\mathbb{R})$ which characterizes completely the elements in $(S)^*$ (as well as the strong convergence in $(S)^*$). In the third section, we will see that the $S$-transform is fundamental in order to define the stochastic integral with respect to the Rosenblatt process.\\
\\
\textbf{Definition 1.1}: Let $\Phi\in (S)^*$. For every function $\xi\in S(\mathbb{R})$, we define the $S$-transform of $\Phi$ by:
$$S(\Phi)(\xi)=\left\langle \left\langle \Phi; :\exp(<;\xi>):\right\rangle\right\rangle$$
where $:\exp(<;\xi>):=\exp(<;\xi>-\frac{||\xi||_{L^2(\mathbb{R})}^2}{2})=\sum_{n=0}^{\infty}\frac{I_n(\xi^{\otimes n})}{n!}\in (S)$.\\
\\
\textbf{Remark 1.2}: For every $\Phi\in (L^2)$, we have:
$$S(\Phi)(\xi)=\mathbb{E}[\Phi :\exp(<;\xi>):]=\mathbb{E}^{\mu_{\xi}}[\Phi]$$
where $\mu_{\xi}$ is the translated infinite dimensional measure defined by:
$$\mu_{\xi}(dx)=\exp(<x;\xi>-\frac{||\xi||_{L^2(\mathbb{R})}^2}{2})\mu(dx).$$
Regarding the $S$-transform, we have the following properties and results:\\
\\
\textbf{Theorem 1.3}:
\begin{enumerate}
\item The $S$-transform is injective. If $\forall \xi\in S(\mathbb{R}), S(\Phi)(\xi)=S(\Psi)(\xi)$ then $\Phi=\Psi$ in $(S)^*$.
\item For $\Phi , \Psi\in (S)^*$ there is a unique element $\Phi\diamond\Psi\in (S)^*$ such that for all $\xi\in S(\mathbb{R}), S(\Psi)(\xi)S(\Phi)(\xi)=S(\Phi\diamond\Psi)(\xi)$. It is called the Wick product of $\Phi$ and $\Psi$.
\item Let $\Phi_n\in (S)^*$ and $F_n=S(\Phi_n)$. Then $\Phi_n$ converges strongly in $(S)^*$ if and only if the following conditions are satisfied:
\begin{itemize}
\item $\underset{n\rightarrow \infty}{\lim}F_n(\xi)$ exists for each $\xi\in S(\mathbb{R})$.
\item There exists strictly positive constants K, a and p independent of n such that:
$$\forall n\in\mathbb{N},\forall\xi\in S(\mathbb{R})\quad|F_n(\xi)|\leq K\exp(a||A^p\xi||_{L^2(\mathbb{R})}^2)$$
\end{itemize}
\end{enumerate}

Next, we introduce the Wick tensors. This object allows us to obtain a certain representation for multiple Wiener-Itô integrals on the white noise space which leads to a formula for the translation of second Wiener chaos elements in any direction $\xi\in S(\mathbb{R})$. For futher details, we refer the reader to chapter 5 of Kuo's book pages 32-35 \cite{MR1387829}.\\
\\
\textbf{Definition 1.4}: The trace operator $\tau$ is the element of $\hat{S}'(\mathbb{R}^2)$ uniquely defined by:
$$\forall\psi,\phi\in S(\mathbb{R})\quad <\tau; \phi\otimes\psi>=<\psi;\phi>$$
\textbf{Definition 1.5}: The Wick tensors of any elements $\omega\in S'(\mathbb{R})$ is defined by:
$$:\omega^{\otimes n}:=\sum_{k=0}^{\lfloor\frac{n}{2}\rfloor}C_{2k}^{n}(2k-1)!!(-1)^k\omega^{\otimes (n-2k)}\hat{\otimes}\tau^{\otimes k}$$
where $C_{2k}^{n}=\frac{n!}{(2k)!(n-2k)!}$, $(2k-1)!!=(2k-1)(2k-3)...3\cdot1$ and $\hat{\otimes}$ is the symmetric tensor product.\\
\\
\textbf{Theorem 1.6}: Let $f\in \hat{L}^2(\mathbb{R}^n)$. $\mu$-almost everywhere, we have:
$$I_n(f)(\omega)=<:\omega^{\otimes n}:;f>$$
\textbf{Proof}: See Theorem 5.4 of Kuo's book \cite{MR1387829}.\\
\\
\textbf{Lemma 1.7}: For all $f\in\hat{L}^2(\mathbb{R}^2)$ and $\xi\in S(\mathbb{R})$, we have:
$$\mu-a.e\quad I_2(f)(\omega+\xi)=I_2(f)(\omega)+2I_1(<f;\xi>)(\omega)+<f;\xi^{\otimes 2}>$$
\textbf{Proof}: For all $\omega \in S'(\mathbb{R})$, we have:
$$:\omega^{\otimes 2}:=\sum_{k=0}^1\dfrac{2}{k!(2-2k)!}\dfrac{(-1)^k}{2^k}\omega^{\otimes (2-2k)}\hat{\otimes}\tau^{\otimes k}$$
$$:\omega^{\otimes 2}:=\omega^{\otimes 2}-\tau$$
where the equalities stand in $S'(\mathbb{R}^2)$. Thus, we have,
$$:(\omega+\xi)^{\otimes 2}:=(\omega+\xi)^{\otimes 2}-\tau$$
$$:(\omega+\xi)^{\otimes 2}:=:\omega^{\otimes 2}:+2\omega\otimes \xi+\xi^{\otimes 2}$$
Let $(f_n)$ be a sequence of functions in $\hat{S}(\mathbb{R}^2)$ converging to $f$ in $\hat{L}^2(\mathbb{R}^2)$. For each $n$, we have:
$$<:(\omega+\xi)^{\otimes 2}:;f_n>=<:\omega^{\otimes 2}:;f_n>+2<\omega\otimes \xi;f_n>+<f_n;\xi^{\otimes 2}>$$
$$<:(\omega+\xi)^{\otimes 2}:;f_n>=<:\omega^{\otimes 2}:;f_n>+2<\omega;<\xi;f_n>>+<f_n;\xi^{\otimes 2}>$$
Thus, passing to the limit, we obtain $\mu$-almost everywhere:
$$I_2(f)(\omega+\xi)=I_2(f)(\omega)+2I_1(<f;\xi>)(\omega)+<f;\xi^{\otimes 2}>.\quad\Box$$

We end this section by introducing the differential operator and its adjoint. These operators play an important role in the context of white noise integration with respect to Brownian motion (see chapter 13.4 of \cite{MR1387829}). Actually, they lead to sufficient condition in order to the white noise integral of a stochastic process to be in $(L)^2$. Thus, as we will see in the third section, they also lead to sufficient condition for the Rosenblatt noise integral of a stochastic process to be a square integrable random variable. We refer the reader to chapter 9 of \cite{MR1387829} for further information concerning these two operators.\\
\\
\textbf{Definition-Theorem 1.8}: Let $y\in S'(\mathbb{R})$ and $\Phi\in(S)$. The operator $D_y$ is continuous from $(S)$ into itself and we have:
$$\forall\omega\in S'(\mathbb{R})\quad D_y(\Phi)(\omega)=\sum_{n=1}^{\infty}n<:\omega:^{\otimes n-1};y\otimes_1\phi_n>,$$
where we denote by $\otimes_1$ the contraction of order $1$ (see \cite{MR2200233}).\\
\\
\textbf{Proof}: See theorem 9.1 of \cite{MR1387829}.$\Box$\\
\\
\textbf{Definition-Theorem 1.9}: Let $y\in S'(\mathbb{R})$ and $\Psi\in(S)^*$. The adjoint operator $D^*_y$ is continuous from $(S)^*$ into itself and we have:
\begin{align*}
\forall\xi\in S(\mathbb{R})\quad S(D^*_y(\Psi))(\xi)=<y;\xi>S(\Psi)(\xi)=S(I_1(y)\diamond\Psi)(\xi)
\end{align*}
where $I_1(y)$ is a generalized Wiener-Itô integral in $(S)^*$. Moreover, we have the following generalized Wiener-Itô decomposition for $D^*_y(\Psi)$:
\begin{align*}
D^*_y(\Psi)(.)=\sum_{n=0}^{\infty}<:.:^{\otimes n+1};y\hat{\otimes}\psi_n>
\end{align*}
\textbf{Proof}: See theorem 9.12, 9.13 and the remark following corollary 9.14 in \cite{MR1387829}.$\Box$

\section{Second order Wiener chaoses and Rosenblatt process.}

In this section, we state the definition of the Rosenblatt process and several known properties of the Rosenblatt distribution (the law of the Rosenblatt process at time $1$). Specifically, we remind the formulae of its characteristic function and of its cumulants. In \cite{Veillette} and more generally in \cite{MR1921743}, it is shown that this characteristic function is infinitely divisible. Moreover in \cite{Maejima20131490}, the authors show that this distribution is selfdecomposable and unimodal. In this section, we show that the translated characteristic function of the ditribution of the Rosenblatt process at $t$ is analytic on $\mathbb{R}$ and admits a certain analytic representation in the neighborhood of the origin which is fundamental in the derivation of the Itô formula for functionals of the Rosenblatt process.\\ 
\\
\textbf{Definition 2.1}: The Rosenblatt process is defined by the following double Wiener-Itô integral representation:
\begin{equation}\forall t\in\mathbb{R}_+\quad X^{H}_t=c(H)\int_{\mathbb{R}^2}\left(\int_{0}^{t}(s-x_1)^{\frac{H}{2}-1}_{+}(s-x_2)^{\frac{H}{2}-1}_{+}ds\right)dB_{x_1}dB_{x_2}
\end{equation}
where $c(H)$ is a normalizing constant such that $\mathbb{E}[|X^H_1|^2]=1$. In particular, by Remark 2.3 of \cite{MR2959876}, we have:
$$c(H)=\sqrt{\dfrac{H(2H-1)}{2}}\dfrac{1}{\beta(1-H,\frac{H}{2})}.$$
\textbf{Remarks 2.2}:
\begin{itemize}
\item Moreover, using the well-known formulae:
$$\forall a,b>0\quad\beta(a,b)=\dfrac{\Gamma(a)\Gamma(b)}{\Gamma(a+b)}\quad\quad \Gamma(a)\Gamma(1-a)=\dfrac{\pi}{\sin(\pi a)},$$
we obtain:
$$c(H)=\sqrt{\dfrac{H(2H-1)}{2}}\Gamma(H)\dfrac{\sin(\pi H)}{\sin(\frac{\pi}{2}H)}\dfrac{1}{(\Gamma(\frac{H}{2}))^2}$$
And thus,
$$X^{H}_t=d(H)\int_{\mathbb{R}^2}\left(\int_{0}^{t}\dfrac{(s-x_1)^{\frac{H}{2}-1}_{+}}{\Gamma(\frac{H}{2})}\dfrac{(s-x_2)^{\frac{H}{2}-1}_{+}}{\Gamma(\frac{H}{2})}ds\right)dB_{x_1}dB_{x_2}$$
\item As we can see the kernel of the Rosenblatt process verifies the following properties:
\begin{enumerate}
\item $f_t^H(x_1,x_2)=d(H)\int_{0}^{t}\frac{(s-x_1)^{\frac{H}{2}-1}_{+}}{\Gamma(\frac{H}{2})}\frac{(s-x_2)^{\frac{H}{2}-1}_{+}}{\Gamma(\frac{H}{2})}ds$ is in $\hat{L}^2(\mathbb{R}^2)$. 
\item $f^H_t(tx_1,tx_2)= t^{H-1}f^H_1(x_1,x_2)$.
\item $f^H_{t+h}(x_1,x_2)-f^H_{t}(x_1,x_2)=f^H_h(x_1-t,x_2-t)$.
\end{enumerate}
\item It is well known that with each function $f\in \hat{L}^2(\mathbb{R}^2)$, we can associate the self-adjoint Hilbert-Schmidt operator, $T_f$, defined on $L^2(\mathbb{R})$ by:
$$\forall g\in L^2(\mathbb{R})\quad T_f(g)(.)=\int_{\mathbb{R}}f(x_1,.)g(x_1)dx_1$$
For the Rosenblatt process, we denote by $T_1$ the self-adjoint Hilbert-Schmidt positive operator associated with $f^H_1$. Denoting by $\lambda_1\geq\lambda_2\geq...\geq\lambda_n\geq...$ and $\{e_n\}$ the eigenvalues and the associated eigenvectors of $T_1$, we obtain the following representation for the kernel, $f_t^H$, in $\hat{L}^2(\mathbb{R}^2)$:
\begin{equation}
f_t^H(.,.)=\sum_{n=1}^{+\infty}\lambda_n t^{H-1}(e_n\otimes e_n)(\frac{.}{t},\frac{.}{t})
\end{equation}
\item From this formula, we classically obtain a peculiar representation in law (in $(L^2)$ and almost surely as well) for the Rosenblatt process at time $t$:
$$I_2(f^H_t)\overset{(d)}{=}t^{H}\sum_{n=1}^{+\infty}\lambda_n(X_n^2-1)$$
where $(X_n)$ is a sequence of independent standard normal random variables.
\item We end this series of remarks by reminding the formula of the cumulants for any second Wiener chaos element (see \cite{MR2970700}):
$$\kappa_1(I_2(f))=0\quad\quad \forall r\geq 2\quad\kappa_r(I_2(f))=2^{r-1}(r-1)!\sum_{n=1}^{+\infty}\lambda_n^{r}(f)$$
where $(\lambda_n(f))$ are the eigenvalues of $T_f$.
\end{itemize}

In the computation of the translated characteristic function for the Rosenblatt process at time $t$, we need the three following lemmas. The first one is the second part of Lemma 2.1 in \cite{MR2959876}. The proof of the second one is quite similar to the calculation of the Rosenblatt distribution's cumulants. The last one is a simple reformulation of Lemma $1$ using the diagonal representation of the kernel $f^H_t$.\\
\\
\textbf{Lemma 2.3}: For every $s,r>0$:
$$\int_{-\infty}^{s\wedge r}(s-x_1)^{\frac{H}{2}-1}_{+}(r-x_1)_{+}^{\frac{H}{2}-1}dx_1=\beta(1-H,\dfrac{H}{2})|s-r|^{H-1}$$
\textbf{Lemma 2.4}: For every $\xi \in S(\mathbb{R})$, we have:
$$<f_t^H;\xi^{\otimes 2}>=d(H)\int_{0}^t[I_+^{\frac{H}{2}}(\xi)(s)]^2ds$$
and,
$$<\underbrace{(...((f_t^H\otimes_1f_t^H)\otimes_1f_t^H)...\otimes_1f_t^H)}_{k-1 \times\otimes_1};\xi^{\otimes 2}>=d(H)\sqrt{\frac{H(2H-1)}{2}}^{k-1}$$
$$\int_{0}^t\int_{0}^tI_+^{\frac{H}{2}}(\xi)(s)I_+^{\frac{H}{2}}(\xi)(r)K^{k-2}_t(s,r)dsdr$$
where $I_+^{H/2}(\xi)(x)=1/(\Gamma(H/2))\int_{-\infty}^x(x-y)^{H/2-1}_+\xi(y)dy$ is the fractional integral of Weyl's type of order $H/2$ of $\xi$ (see \cite{samko1993fractional} chapter $2$ section $5.1$) and $K^{k-2}_t(s,r)$ is the sequence of kernels defined by:
$$K_t^0(s,r)=|s-r|^{H-1}\quad K_t^1(s,r)=\int_{0}^t|s-u|^{H-1}|r-u|^{H-1}du$$
$$\forall k\geq 3, K^{k-2}_t(s,r)=\int_0^t...\int_0^t|s-x_1|^{H-1}|x_2-x_1|^{H-1}...|x_{k-2}-x_{k-3}|^{H-1}|r-x_{k-2}|^{H-1}dx_1...dx_{k-2}$$
\textbf{Remark 2.5}: The type of integrals used in the definition of the sequence of kernels $(K^k_t(s,r))_{k}$ intervenes in the calculation of the cumulants of the Rosenblatt distribution (see \cite{Veillette}). Indeed, we have, for all $r\geq 2$:
\begin{equation}
\kappa_r(X^H_1)=2^{r-1}(r-1)!\left(\sqrt{\frac{H(2H-1)}{2}}\right)^{r}\int_{0}^1...\int_{0}^1|x_1-x_2|^{H-1}...|x_r-x_{r-1}|^{H-1}|x_1-x_r|^{H-1}dx_1...dx_r
\end{equation}
And, by Remark 2.2,
\begin{equation}\sum_{n=1}^{+\infty}\lambda_n^r=\left(\sqrt{\frac{H(2H-1)}{2}}\right)^{r}\int_{0}^1...\int_{0}^1|x_1-x_2|^{H-1}...|x_r-x_{r-1}|^{H-1}|x_1-x_r|^{H-1}dx_1...dx_r
\end{equation}
\textbf{Proof}:\\ 
Let $\xi$ be in $S(\mathbb{R})$. Note that $I_+^{\frac{H}{2}}(\xi)\in C^{\infty}(\mathbb{R})$. Using Fubini's theorem and the definition of the fractional integral of Weyl's type of order $\frac{H}{2}$, we have:
$$<f^H_t;\xi^{\otimes 2}>=d(H)\int_{0}^{t}[I_+^{\frac{H}{2}}(\xi)(s)]^2ds$$
Moreover, using Fubini's theorem and Lemma $2.3$, we have:
$$<f^H_t\otimes_1 f^H_t;\xi^{\otimes 2}>=\int_{\mathbb{R}^2}(f_t^H\otimes_1 f_t^H)(x_1,x_2)\xi(x_1)\xi(x_2)dx_1dx_2$$
$$=\int_{\mathbb{R}^2}\left(\int_{\mathbb{R}}f^H_t(x_1,u)f^H_t(x_2,u)du\right)\xi(x_1)\xi(x_2)dx_1dx_2$$
$$=(d(H))^2\int_0^t\int_0^tI_+^{\frac{H}{2}}(\xi)(s)I_+^{\frac{H}{2}}(\xi)(r)\left(\int_{\mathbb{R}}\dfrac{(s-u)^{\frac{H}{2}-1}_+}{\Gamma(\frac{H}{2})}\dfrac{(r-u)^{\frac{H}{2}-1}_+}{\Gamma(\frac{H}{2})}du\right)dsdr$$
$$=(d(H))^2\dfrac{\beta(1-H,\frac{H}{2})}{(\Gamma(\frac{H}{2}))^2}\int_0^t\int_0^tI_+^{\frac{H}{2}}(\xi)(s)I_+^{\frac{H}{2}}(\xi)(r)|s-r|^{H-1}dsdr$$
$$=d(H)\sqrt{\dfrac{H(2H-1)}{2}}\int_0^t\int_0^tI_+^{\frac{H}{2}}(\xi)(s)I_+^{\frac{H}{2}}(\xi)(r)|s-r|^{H-1}dsdr$$
Now, using iteratively Fubini's theorem and Lemma $2.3$, we obtain straightforwardly for all $k\geq 3$:
$$<\underbrace{(...((f_t^H\otimes_1f_t^H)\otimes_1f_t^H)...\otimes_1f_t^H)}_{k-1 \times\otimes_1};\xi^{\otimes 2}>=d(H)\left(\sqrt{\frac{H(2H-1)}{2}}\right)^{k-1}$$
$$\int_{0}^t\int_{0}^tI_+^{\frac{H}{2}}(\xi)(s)I_+^{\frac{H}{2}}(\xi)(r)K^{k-2}_t(s,r)dsdr.\quad\Box$$
\\
\textbf{Lemma 2.6}: For any $\xi\in S(\mathbb{R})$: 
$$I_2(f^H_t)(.+\xi)\overset{(d)}{=}t^H\sum_{n=1}^{+\infty}\lambda_n\left[X^2_n-1+2<\xi;\dfrac{e_n(./t)}{\sqrt{t}}>X_n\right]+<f^H_t;\xi^{\otimes 2}>$$
where $(X_n)$ is a sequence of independent standard normal random variables.\\
\\
\textbf{Proof}: This is just a direct application of Lemma 1.7 and of the diagonal decomposition of the kernel $f^H_t$. $\Box$\\
\\
\textbf{Definition 2.7}: A Rosenblatt distribution is the law of the Rosenblatt process at time $1$ whose characteristic function admits the following analytic representation in a neighborhood of the origin:
\begin{equation}
\quad\mathbb{E}[\exp(i\theta X^H_1)]=\exp\left(\frac{1}{2}\sum_{k=2}^{+\infty}\frac{(2i\theta)^k}{k}\left(\sqrt{\frac{(2H-1)H}{2}}\right)^kc_k\right)
\end{equation}
where,
$$c_k=\int_{0}^1...\int_{0}^1|x_1-x_2|^{H-1}...|x_k-x_{k-1}|^{H-1}|x_1-x_k|^{H-1}dx_1...dx_k$$

\textbf{Remark 2.8}:
\begin{itemize}
\item Using Cauchy-Schwarz inequality, one can show that $|c_k|\leq (\frac{1}{H(2H-1)})^{\frac{k}{2}}$ ensuring that the former series is convergent in a neighborhood of the origin and thus analytic in this neighborhood. Nevertheless, using a result by Lukacs (see theorem 7.1.1 in \cite{MR0346874}), this characteristic function is actually analytic in a horizontal strip of the complex plane containing the real line.
\item In \cite{Veillette} and \cite{MR1921743} this distribution is shown to be infinitely divisible and, thus, its characteristic function admits a Lévy-Khintchine representation. Using this representation, one can show that the density of the Rosenblatt distribution exists and is infinitely differentiable with all derivative tending to $0$ at $+\infty$ (Corollary 4.3 of \cite{Veillette}).
\end{itemize}
We are now ready to state the main result of this section.\\
\\
\textbf{Theorem 2.9}: Let $\xi\in S(\mathbb{R})$ and $t>0$. We have (for $\theta$ being small enough):
$$
\mathbb{E}^{\mu_{\xi}}[\exp(i\theta X^H_t)]=\exp(i\theta<f^H_t;\xi^{\otimes 2}>)\exp\left(\frac{1}{2}\sum_{k=2}^{+\infty}\frac{(2i\theta t^H)^k}{k}\left(\sqrt{\frac{(2H-1)H}{2}}\right)^kc_k\right)$$
$$\times\exp\left(\sum_{k=2}^{+\infty}(2i\theta)^k||T^{\frac{k}{2}}_t(\xi)||_{L^2(\mathbb{R})}^2\right)$$
where $T^{k/2}_t$ denotes the operator $T_t$, associated with $f^H_t$, iterated $k/2$-times. Moreover, we have:
\begin{equation}
\forall k\geq 2\quad||T^{\frac{k}{2}}_t(\xi)||_{L^2(\mathbb{R})}^2=\sum_{n=1}^{\infty}(\lambda_nt^{H})^k(<\xi;\frac{e_n(./t)}{\sqrt{t}}>)^2
\end{equation}
$$=<\underbrace{(...((f_t^H\otimes_1f_t^H)\otimes_1f_t^H)...\otimes_1f_t^H)}_{k-1 \times\otimes_1};\xi^{\otimes 2}>$$
\textbf{Remark 2.10}:
\begin{itemize}
\item Under the probability $\mu_{\xi}$, the random variable $X^H_t$ has a non-zero mean and is the sum of two independent random variables with one of them having the law of the Rosenblatt process at time $t$.
\item Trivially, we obtain the following estimates ensuring the convergence around the origin of the second series in the expression of $\mathbb{E}^{\mu_{\xi}}[\exp(i\theta X^H_t)]$:
$$||T^{\frac{k}{2}}_t(\xi)||_{L^2(\mathbb{R})}^2\leq ||\underbrace{(...((f_t^H\otimes_1f_t^H)\otimes_1f_t^H)...\otimes_1f_t^H)}_{k-1 \times\otimes_1}||_{L^2(\mathbb{R}^2)}||\xi||_{L^2(\mathbb{R})}^2$$
$$||T^{\frac{k}{2}}_t(\xi)||_{L^2(\mathbb{R})}^2\leq \left(\sum_{n=1}^{\infty}(\lambda_n t^H)^{2k}\right)^{\frac{1}{2}}||\xi||_{L^2(\mathbb{R})}^2=t^{Hk}\left(\sqrt{\frac{(2H-1)H}{2}}\right)^{k}(c_{2k})^{\frac{1}{2}}||\xi||_{L^2(\mathbb{R})}^2$$
\end{itemize}
\textbf{Proof}: Let $\xi\in S(\mathbb{R})$ and $t>0$. By Lemma $2.6$, we have:
$$X^H_t(.+\xi)\overset{(d)}{=}t^H\sum_{n=1}^{+\infty}\lambda_n\left[X^2_n-1+2<\xi;\dfrac{e_n(./t)}{\sqrt{t}}>X_n\right]+<f^H_t;\xi^{\otimes 2}>$$
Let us compute the characteristic function of $X^2_n-1+2<\xi;\dfrac{e_n(./t)}{\sqrt{t}}>X_n$ where $X_n$ is a standard normal random variable. Let $\alpha>0$, $\beta\in\mathbb{R}$ and $X$ be a centered normal random variable with variance $\frac{1}{2\alpha}$. By standard computations, we get the following formulae for the density and the Laplace transform of $Y=(X+\beta)^2$: 
$$\forall y\in\mathbb{R}\quad f_Y(y)=1_{(0,+\infty)}(y)\sqrt{\dfrac{\alpha}{\pi}}\exp(-\alpha\beta^2)\dfrac{\exp(-\alpha y)}{\sqrt{y}}\cosh(2\alpha\beta\sqrt{y})$$
$$\forall\mu>0\quad \mathbb{E}[\exp(-\mu Y)]=\left(\dfrac{\alpha}{\alpha+\mu}\right)^{\frac{1}{2}}\exp(-\alpha\beta^2)\exp(\frac{\alpha^2\beta^2}{\alpha+\mu})$$
Subtracting $(\frac{1}{2\alpha}+\beta^2)$ from $Y$ in order to obtain a centered random variable, we have:
$$\forall\mu>0\quad \mathbb{E}[\exp(-\mu Y')]=\left(\dfrac{\alpha}{\alpha+\mu}\right)^{\frac{1}{2}}\exp(\frac{\mu}{2\alpha})\exp(-\alpha\beta^2)\exp(\frac{\alpha^2\beta^2}{\alpha+\mu})\exp(\mu\beta^2)$$
Putting $\alpha=\frac{1}{2}$, $\beta_n=<\xi;\frac{e_n(./t)}{\sqrt{t}}>$ and $Y=(X_n+\beta_n)^2$, we get:
$$\forall n\in\mathbb{N}^*,\forall\mu>0\quad\mathbb{E}[\exp(-\mu\lambda_n(X_n^2-1+2\beta_nX_n))]=\left(\dfrac{1}{1+2\lambda_n\mu}\right)^{\frac{1}{2}}\exp(\lambda_n\mu)\exp(\frac{2\beta_n^2\lambda_n^2\mu^2}{1+2\lambda_n\mu})$$
Leading to the following formula for the characteristic function of $\lambda_n(X^2_n-1+2\beta_nX_n)$:
$$\forall\theta\in\mathbb{R}\quad\mathbb{E}[\exp(i\theta\lambda_n(X_n^2-1+2\beta_nX_n))]=\left(\dfrac{1}{1-2i\lambda_n\theta}\right)^{\frac{1}{2}}\exp(-i\lambda_n\theta)\exp\left(-\frac{2\beta_n^2\lambda_n^2\theta^2}{1-2i\lambda_n\theta}\right)$$
Using the independence of the $X_n$'s, we have for any $N\geq 1$:
$$\mathbb{E}[\exp(i\theta t^H\sum_{n=1}^N\lambda_n(X_n^2-1+2\beta_nX_n))]=\left(\prod_{n=1}^N\left(\dfrac{1}{1-2i\lambda_nt^H\theta}\right)^{\frac{1}{2}}\exp(-i\lambda_nt^H\theta)\right)\exp\left(-2\theta^2\sum_{n=1}^N\frac{\beta_n^2\lambda_n^2t^{2H}}{1-2i\lambda_nt^H\theta}\right)$$
The first product term converges to the distribution of the Rosenblatt process at time $t$ as $N\rightarrow+\infty$. Using the following estimate:
$$\forall\theta\in\mathbb{R}\quad|\dfrac{\beta_n^2\lambda_n^2}{1-2i\lambda_nt^H\theta}|\leq |\beta_n|^2\lambda_n^2$$
and, (6) for $k=2$, we obtain:
$$\forall\theta\in\mathbb{R}\quad\mathbb{E}^{\mu_{\xi}}[\exp(i\theta X^H_t)]=\exp(i\theta<f^H_t;\xi^{\otimes 2}>)\left(\prod_{n=1}^{+\infty}\left(\dfrac{1}{1-2i\lambda_nt^H\theta}\right)^{\frac{1}{2}}\exp(-i\lambda_nt^H\theta)\right)$$
$$\times\exp\left(-2\theta^2\sum_{n=1}^{+\infty}\frac{\beta_n^2\lambda_n^2t^{2H}}{1-2i\lambda_nt^H\theta}\right)$$
For $\theta$ small enough, using (4) and (6), we get:
$$\prod_{n=1}^{+\infty}\left(\left(\dfrac{1}{1-2i\lambda_nt^H\theta}\right)^{\frac{1}{2}}\exp(-i\lambda_nt^H\theta)\right)=\exp\left(\frac{1}{2}\sum_{k=2}^{+\infty}\frac{(2i\theta t^H)^k}{k}\left(\sqrt{\frac{(2H-1)H}{2}}\right)^kc_k\right)$$
$$\exp\left(-2\theta^2\sum_{n=1}^{+\infty}\frac{\beta_n^2\lambda_n^2t^{2H}}{1-2i\lambda_nt^H\theta}\right)=\exp\left(\sum_{k=2}^{+\infty}(2i\theta)^k||T^{\frac{k}{2}}_t(\xi)||^2\right)$$
Ensuring the analyticity of the translated characteristic function in a neighborhood of the origin, and, by theorem 7.1.1 of \cite{MR0346874}, the analyticity in a certain horizontal strip. Thus, this characteristic function is in $C^{\infty}(\mathbb{R})$.\\
\\ 
Moreover, using (2) and (6), we have:
$$||T^{\frac{k}{2}}_t(\xi)||_{L^2(\mathbb{R})}^2=<\underbrace{(...((f_t^H\otimes_1f_t^H)\otimes_1f_t^H)...\otimes_1f_t^H)}_{k-1 \times\otimes_1};\xi^{\otimes 2}>.\Box$$

\section{Stochastic calculus and Itô formula with respect to the Rosenblatt process.}

In this section, we define differentiability as well as integrability of stochastic distribution processes in $(S)^*$. In particular, we compute the $(S)^*$ derivative of the Rosenblatt process, $\dot{X}^H_t$, named the Rosenblatt noise. Thus, using Wick product and Pettis integration, we define the Rosenblatt noise integral of integrable stochastic distribution process. Moreover, we derive a sufficient condition for the Rosenblatt noise integral of a stochastic process to be an element in $(L^2)$. We end this section with Itô formulae for functionals of the Rosenblatt process.\\
\\
\textbf{Definition 3.1}:
\begin{enumerate}
\item Let $I\subset\mathbb{R}$ be an interval. A mapping $X:I\rightarrow (S)^*$ is called a stochastic distribution process.
\item A stochastic distribution process $X$ is said to be differentiable, if $\underset{h\rightarrow 0}{\lim}\frac{X_{t+h}-X_t}{h}$ exists in $(S)^*$.
\end{enumerate}
\textbf{Remark 3.2}: By the last point of Theorem 1.3, convergence in $(S)^*$ is ensured by pointwise convergence of the $S$-transform and a uniform growth condition.\\
\\
\textbf{Lemma 3.3}: Let $\xi\in S(\mathbb{R})$ and $H\in(\frac{1}{2},1)$. Then $I^{\frac{H}{2}}_+(\xi)\in C^{\infty}(\mathbb{R})$ and it exists a strictly positive constant $C_H$ such that:
$$||I^{\frac{H}{2}}_+(\xi)||_{\infty}\leq C_H(||\xi||_{\infty}+||\xi^{(1)}||_{\infty}+||\xi||_1)$$
\textbf{Proof}: See the proof of Theorem 2.3 in \cite{MR1956473}.$\Box$\\
\\
\textbf{Lemma 3.4}: The Rosenblatt process is $(S)^*$ differentiable and its derivative, the Rosenblatt noise, admits the following generalized double Wiener-Itô integral representation (as an element of $(S)^*$):
$$\forall t>0\quad \dot{X}^H_t=d(H)I_2(\delta_t^{\otimes 2}\circ (I^{\frac{H}{2}}_+)^{\otimes 2})$$
And its $S$-transform is equal to:
$$\forall\xi\in S(\mathbb{R})\quad S(\dot{X}^H_t)(\xi)=d(H)(I^{\frac{H}{2}}_+(\xi)(t))^2$$
\textbf{Proof}: Let $t>0$. By Fubini's theorem, we have:
$$\forall \xi\in S(\mathbb{R})\quad S(X^H_t)(\xi)=d(H)\int_0^t(I^{\frac{H}{2}}_+(\xi)(s))^2ds$$
Thus,
$$\forall h>0\quad S(\dfrac{X^H_{t+h}-X^H_{t}}{h})(\xi)=\frac{d(H)}{h}\int_t^{t+h}(I^{\frac{H}{2}}_+(\xi)(s))^2ds$$
By Lemma 3.3, we obtain the following estimate:
$$|S(\dfrac{X^H_{t+h}-X^H_{t}}{h})(\xi)|\leq C_Hd(H)(||\xi||_{\infty}+||\xi^{(1)}||_{\infty}+||\xi||_1)^2$$
Using Theorem 1.3, $\frac{X^H_{t+h}-X^H_{t}}{h}$ converges in $(S)^*$ and the limit is completely characterized by its $S$-transform:
$$S(\dot{X}^H_t)(\xi)=d(H)(I^{\frac{H}{2}}_+(\xi)(t))^2$$
By a similar argument than the one used in the representation of the fractional noise (see \cite{MR1956473} page $90$), we can extend the double Wiener-Itô integral to elements of $\hat{S}'(\mathbb{R}^2)$. Consequently, we get the following representation for the Rosenblatt noise:
$$\dot{X}^H_t=d(H)I_2(\delta_t^{\otimes 2}\circ (I^{\frac{H}{2}}_+)^{\otimes 2})$$
\textbf{Remark 3.5}: Let $d\geq 1$ and let $H\in (\frac{1}{2},1)$. Let $H_0=\frac{1}{2}+\frac{H-1}{d}$. We define the Hermite process of order $d$ by the following multiple Wiener-Itô integral representation:
$$\forall t>0\quad Y^{H,d}_t=d(H_0)\int_{\mathbb{R}}...\int_{\mathbb{R}}\left(\int_0^t\prod_{j=1}^d\left(\dfrac{(s-x_j)^{H_0-1}_+}{\Gamma(H_0)}\right)ds\right)dB_{x_1}...dB_{x_d},$$
where $d(H_0)$ is a normalizing constant such that $\mathbb{E}[|Y^{H,d}_1|^2]=1$. The Hermite processes of order $1$ and $2$ are respectively fractional Brownian motion and the Rosenblatt process. As previously, we obtain the following formula and representation for the derivative in $(S)^*$ of the Hermite process of order $d$:
$$\forall t>0\quad \dot{Y}^{H,d}_t=d(H_0)I_d(\delta_t^{\otimes d}\circ(I^{H_0}_+)^{\otimes d}),$$
$$\forall\xi\in S(\mathbb{R})\quad S(\dot{Y}^{H,d}_t)(\xi)=d(H_0)(I^{H_0}_+(\xi)(t))^d.$$
\textbf{Definition 3.6}: Let $k\geq 2$. For any $t>0$, we define the following sequence of stochastic processes (belonging to the second Wiener chaos):
$$X^{H,k}_t=\int_{\mathbb{R}}\int_ {\mathbb{R}}\underbrace{(...((f_t^H\otimes_1f_t^H)\otimes_1f_t^H)...\otimes_1f_t^H)}_{k-1 \times\otimes_1}(x_1,x_2)dB_{x_1}dB_{x_2}$$
Moreover, the $S$-transform of this process is given by:
$$\forall\xi\in S(\mathbb{R})\quad S(X^{H,k}_t)(\xi)=d(H)\sqrt{\frac{H(2H-1)}{2}}^{k-1}
\int_{0}^t\int_{0}^tI_+^{\frac{H}{2}}(\xi)(s)I_+^{\frac{H}{2}}(\xi)(r)K^{k-2}_t(s,r)dsdr$$
\textbf{Lemma 3.7}: For any $k\geq 2$, the process $\{X^{H,k}_t:t>0\}$ is $(S)^*$ differentiable and its derivative, $\dot{X}^{H,k}_t$, is uniquely defined by the following $S$-transform:
\begin{align*}
\forall\xi\in S(\mathbb{R})\quad S(\dot{X}^{H,k}_t)(\xi)=&d(H)\sqrt{\frac{H(2H-1)}{2}}^{k-1}\\
\Big[&(H(k-1)+1)t^{H(k-1)}\int_0^1\int_0^1I_+^{\frac{H}{2}}(\xi)(tu)I_+^{\frac{H}{2}}(\xi)(tv)K^{k-2}_1(u,v)dudv
\\&+t^{H(k-1)+1}\int_0^1\int_0^1\dfrac{d}{dt}[I_+^{\frac{H}{2}}(\xi)(tu)I_+^{\frac{H}{2}}(\xi)(tv)]K^{k-2}_1(u,v)dudv\Big].
\end{align*}
\textbf{Proof}: Let $(a,b)\in \mathbb{R}^*_+$ such that $0<a<b<\infty$. Let $t\in (a,b)$ and $\xi\in S(\mathbb{R})$. In order to prove the $(S)^*$ differentiability of the process $X^{H,k}_t$, we will proceed as in the proof of the former lemma. Let $h>0$ such that $t+h\in (a,b)$. We have:
$$S(X^{H,k}_t)(\xi)=d(H)\sqrt{\frac{H(2H-1)}{2}}^{k-1}
\int_{0}^t\int_{0}^tI_+^{\frac{H}{2}}(\xi)(s)I_+^{\frac{H}{2}}(\xi)(r)K^{k-2}_t(s,r)dsdr$$
Using scaling property of the sequence of kernels $K^{k-2}_t(s,r)$, we get:
$$S(X^{H,k}_t)(\xi)=d(H)\sqrt{\frac{H(2H-1)}{2}}^{k-1}t^{H(k-1)+1}\int_0^1\int_0^1I_+^{\frac{H}{2}}(\xi)(tu)I_+^{\frac{H}{2}}(\xi)(tv)K^{k-2}_1(u,v)dudv$$
Thus, we have:
$$|S(X^{H,k}_{t+h})(\xi)-S(X^{H,k}_t)(\xi)|\leq d(H)\sqrt{\frac{H(2H-1)}{2}}^{k-1}\left[(I)+(II)\right]$$
where,
$$(I)=\left|(t+h)^{H(k-1)+1}-t^{H(k-1)+1}\right|\left|\int_0^1\int_0^1I_+^{\frac{H}{2}}(\xi)((t+h)u)I_+^{\frac{H}{2}}(\xi)((t+h)v)K^{k-2}_1(u,v)dudv\right|$$
$$(II)=t^{H(k-1)+1}\left|\int_0^1\int_0^1\left[I_+^{\frac{H}{2}}(\xi)((t+h)u)I_+^{\frac{H}{2}}(\xi)((t+h)v)-I_+^{\frac{H}{2}}(\xi)(tu)I_+^{\frac{H}{2}}(\xi)(tv)\right]K^{k-2}_1(u,v)dudv\right|$$
Combining Lemma 3.3 and the Mean value theorem, we obtain for $(I)$:
\begin{equation}
(I)\leq (H(k-1)+1)b^{H(k-1)}C_H^2(||\xi||_{\infty}+||\xi^{(1)}||_{\infty}+||\xi||_1)^2\underbrace{\left(\int_0^1\int_0^1K^{k-2}_1(u,v)dudv\right)}_{<\infty}h
\end{equation}
For the second term, we have:
$$(II)\leq 2b^{H(k-1)+1}||I^{\frac{H}{2}}_+(\xi)||_{\infty}\int_0^1\int_0^1\left|I_+^{\frac{H}{2}}(\xi)((t+h)u)-I_+^{\frac{H}{2}}(\xi)(tu)\right|K^{k-2}_1(u,v)dudv$$
Using the fact that $I^{\frac{H}{2}}_+(\xi)\in C^{\infty}(\mathbb{R})$, we have as previously:
$$(II)\leq 2b^{H(k-1)+1}||I^{\frac{H}{2}}_+(\xi)||_{\infty}||I^{\frac{H}{2}}_+(\xi^{(1)})||_{\infty}\left(\int_0^1\int_0^1K^{k-2}_1(u,v)dudv\right)h$$
\begin{align}
(II)\leq 2b^{H(k-1)+1}&C_H^2(||\xi||_{\infty}+||\xi^{(1)}||_{\infty}+||\xi||_1)\nonumber\\
&(||\xi^{(1)}||_{\infty}+||\xi^{(2)}||_{\infty}+||\xi^{(1)}||_1)\left(\int_0^1\int_0^1K^{k-2}_1(u,v)dudv\right)h
\end{align}
By $(7)$ and $(8)$, $|S(\frac{X^{H,k}_{t+h}-X^{H,k}_t}{h})(\xi)|$ verifies the uniform growth condition of the last point of Theorem 1.3.\\
Moreover, the pointwise convergence of $S(\frac{X^{H,k}_{t+h}-X^{H,k}_t}{h})(\xi)$ as $h\rightarrow 0^+$ is ensured by the dominated convergence theorem.\\ 
Consequently, the process $X^{H,k}$ is $(S)^*$ differentiable for any $k\geq 2$ and its $(S)^*$-derivative admits the following $S$-transform:
\begin{align*}
S(\dot{X}^{H,k}_t)(\xi)=d(H)\sqrt{\frac{H(2H-1)}{2}}^{k-1}\Big[&(H(k-1)+1)t^{H(k-1)}\int_0^1\int_0^1I_+^{\frac{H}{2}}(\xi)(tu)I_+^{\frac{H}{2}}(\xi)(tv)K^{k-2}_1(u,v)dudv\\
&+t^{H(k-1)+1}\int_0^1\int_0^1\dfrac{d}{dt}[I_+^{\frac{H}{2}}(\xi)(tu)I_+^{\frac{H}{2}}(\xi)(tv)]K^{k-2}_1(u,v)dudv\Big].\Box
\end{align*}

As well as differentiability, we can define integrability of a stochastic distribution process in term of the $S$-transform.\\
\\
\textbf{Definition 3.8}: A stochastic distribution process $X:I\rightarrow (S)^*$ is integrable if:
\begin{enumerate}
\item $\forall\xi\in S(\mathbb{R})$, $S(X_.)(\xi)$ is measurable on $I$.
\item $\forall\xi\in S(\mathbb{R})$, $S(X_.)(\xi)\in L^1(I)$.
\item $\int_IS(X_t)(\xi)dt$ is the $S$-transform of a certain Hida distribution. 
\end{enumerate}
We have the following criterion for integrability of stochastic distribution process.\\
\\
\textbf{Theorem 3.9}: Let $X:I\rightarrow (S)^*$ be a stochastic distribution process satisfying:
\begin{enumerate}
\item $\forall\xi\in S(\mathbb{R})$, $S(X_.)(\xi)$ is measurable on $I$.
\item There is a $p\in\mathbb{N}$, a strictly positive constant $a$ and a non-negative function $L\in L^1(I)$ such that:
$$\forall\xi\in S(\mathbb{R})\quad |S(X_t)(\xi)|\leq L(t)\exp\left(a||A^p\xi||^2_2\right)$$
\end{enumerate}
Then $X$ is $(S)^*$-integrable.\\
\\
\textbf{Proof}: See theorem 13.5 in \cite{MR1387829}.$\Box$\\
\\
We are now ready to state the definition of the Rosenblatt noise integral. This definition is motivated by the one of the white noise integral. Indeed, in the case of stochastic integration with respect to Brownian motion in a white noise context, this integral coincides with the classical Itô integral for $L^2(\Omega\times I)$ nonanticipating integrands (see theorem 13.12 in \cite{MR1387829}).\\
\\
\textbf{Definition-Theorem 3.10}: Let $\{\phi_t;t\in I\}$ be a $(S)^*$ stochastic process satisfying the assumptions of Theorem 3.9. Then, $\phi_t\diamond \dot{X}^H_t$ is $(S)^*$ integrable over $I$ and we define the Rosenblatt noise integral of $\{\phi_t\}$ by:
\begin{align*}
\int_I\phi_tdX^H_t=\int_I\phi_t\diamond \dot{X}^H_tdt
\end{align*}
Moreover, we have the following representation:
\begin{align*}
\int_I\phi_tdX^H_t=\int_I (D^*_{\sqrt{d(H)}\delta_t\circ I^{\frac{H}{2}}_+})^2(\phi_t)dt
\end{align*}
\textbf{Proof}: Using Lemma 3.4 and Theorem 3.9, one can show that $\phi_t\diamond \dot{X}^H_t$ is $(S)^*$-integrable over $I$. Let us compute the $S$-transform of $\int_I\phi_tdX^H_t$.
\begin{align*}
\forall\xi\in S(\mathbb{R})\quad S(\int_I\phi_tdX^H_t)(\xi)&=\int_IS(\phi_t)(\xi)S(\dot{X}^H_t)(\xi)dt\\
&=\int_IS(\phi_t)(\xi)d(H)(I^{\frac{H}{2}}_+(\xi)(t))^2dt
\end{align*}
Then, using Definition-Theorem 1.9 iteratively, we obtain:
\begin{align*}
S(\int_I\phi_tdX^H_t)(\xi)=\int_IS((D^*_{\sqrt{d(H)}\delta_t\circ I^{\frac{H}{2}}_+})^2(\phi_t))(\xi)dt.\quad\Box
\end{align*}
\\
Before stating the next theorem, we need to introduce some notations. First of all, we need to consider the following multidimensional version of fractional integrals (chapter $5$, paragraph $24$ in \cite{samko1993fractional}):
\begin{align*}
I^{\alpha_1,...,\alpha_m}_{+,...,+}(\phi)=(I^{\alpha_1}_+\otimes...\otimes I^{\alpha_m}_+)(\phi)
\end{align*}
where $m\geq 1$, $\alpha_i\in(0,1)$ for every $i\in\{1,...,m\}$ and $\phi\in S(\mathbb{R}^m)$. By Theorem $24.1$ of \cite{samko1993fractional}, we know that the operator $I^{\alpha_1,...,\alpha_m}_{+,...,+}$ is bounded from $L^{(p_1,...,p_m)}(\mathbb{R}^m)$ to $L^{(q_1,...,q_m)}(\mathbb{R}^m)$ if and only if, $1<p_i<1/\alpha_i$ and $q_i=p_i/(1-\alpha_ip_i)$, for every $i\in\{1,...,m\}$. $L^{(p_1,...,p_m)}(\mathbb{R}^m)$ refers to the mixed-norm space of index $(p_1,...,p_m)$ (see chapter $2.48$ in \cite{MR2424078}). When some of the $\alpha_i$ are equal to $0$, we replace $I^{\alpha_i}_+$ by the identity operator and the theorem is still true with $q_i=p_i$ and $p_i\geq 1$. In the sequel, for fixed $m\geq 2$ and $H\in(1/2,1)$, we will have to consider the operators $I^{(0,...,0,H/2)}_{+,...,+}$ and $I^{(0,...,0,H/2,H/2)}_{+,...,+}$ from $L^2(\mathbb{R}^m)$ to $L^{(2,...,2,2/(1-H))}(\mathbb{R}^m)$ and $L^{(2,...,2,2/(1-H),2/(1-H))}(\mathbb{R}^m)$ respectively. Moreover, we will have to consider the restriction of $I^{(0,...,0,H/2,H/2)}_{+,...,+}(f)$ to the hyperplane $\Delta_m=\{(x_1,...,x_m)\in\mathbb{R}^m:\ x_{m-1}=x_m\}$ for functions $f$ not smooth enough. For this purpose, we introduce the function spaces with dominating mixed smoothness which were defined by Lizorkin and Nikol'skii in \cite{MR0192223} (see \cite{MR2657119} for a recent survey on these spaces and their generalizations). For any $(\alpha_1,...,\alpha_m)$ such that $\alpha_i\geq 0$ for every $i\in\{1,...,m\}$, we define $S^{(\alpha_1,...,\alpha_m)}_2H(\mathbb{R}^m)$ by:
\begin{align*}
S^{(\alpha_1,...,\alpha_m)}_2H(\mathbb{R}^m)=\{f\in S'(\mathbb{R}^m): \mathcal{F}^{-1}(\prod_{i=1}^{m}(1+|\xi_i|^2)^{\frac{\alpha_i}{2}}\mathcal{F}(f))\in L^2(\mathbb{R}^m)\}.
\end{align*}
These spaces admit several representations, one of them being particularly relevant in our context (see \cite{MR0192223}). Indeed, using partial and mixed Liouville derivatives of fractional orders (see chapter $5$, paragraph $24.2$ in \cite{samko1993fractional}), we have:
\begin{align*}
S^{(\alpha_1,...,\alpha_m)}_2H(\mathbb{R}^m)=\{f\in L^2(\mathbb{R}^m): \forall \beta=(\beta_1,...,\beta_d), \beta_i=0\ or\ \alpha_i,\  D_{+,...,+}^{(\beta_1,...,\beta_m)}(f)\in L^2(\mathbb{R}^m)\}.
\end{align*}
The mixed Liouville derivative of order $(\alpha_1,...,\alpha_m)$ plays here a dominant role. The trace problem for function spaces with dominating mixed smoothness appears to be much more difficult than for standard isotropic fractional Sobolev spaces. Indeed, function spaces with dominating mixed smoothness are not rotationally invariant in general. For example, trace problems for hyperplane parallel to coordinate axis are different of trace problems for hyperplane in oblique position. Nevertheless, for $m=2$, regarding trace problem on the diagonal of $\mathbb{R}^2$, we have the following result (see Theorem $5.2$ in section 5.1 of \cite{MR2657119} and also \cite{MR2274153} for the proof):\\
\\
\textbf{Theorem 3.11}: Let $\alpha_1$, $\alpha_2$ be strictly positive such that $\alpha_i<1/2$ for every $i=1,2$ and $\alpha_1+\alpha_2>1/2$. Then, we have:
\begin{align*}
Tr \big(S^{(\alpha_1,\alpha_2)}_2H(\mathbb{R}^2)\big)=H^{\alpha_1+\alpha_2-\frac{1}{2}}(\mathbb{R}),
\end{align*}
where $Tr$ is the trace operator on the diagonal of $\mathbb{R}^2$ and $H^{\alpha_1+\alpha_2-1/2}(\mathbb{R})$ is the standard isotropic fractional Sobolev spaces of order $\alpha_1+\alpha_2-1/2$, $W^{\alpha_1+\alpha_2-1/2, 2}(\mathbb{R})$ (see \cite{MR2424078} chapter $7$).\\
\\
In the sequel, we are interested in functions of the form $g=I^{(0,...,0,H/2,H/2)}_{+,...,+}(f)$ for some $f\in L^2(\mathbb{R}^m)$. However, such a function $g$ does not belong to $S^{(0,...,0,H/2,H/2)}_2H(\mathbb{R}^m)$ even if $D^{(0,...,0,H/2,H/2)}_{+...+}(g)\in L^2(\mathbb{R}^m)$ (see Theorem $24.3$ in \cite{samko1993fractional}). This is due to the fact that $I^{(0,...,0,H/2,H/2)}_{+,...,+}(.)$ maps continuously $L^2(\mathbb{R}^m)$ into $L^{(2,...,2,2/(1-H),2/(1-H))}(\mathbb{R}^m)$. Thus, we make the following assumption for the Chaos kernels of the integrand stochastic process which is considered in Theorem 3.13:
\begin{align*}
\mathbf{H1}:\ \forall m\geq 2,\ \forall t\in(a,b),\ I^{(0,...,0,\frac{H}{2},\frac{H}{2})}_{+,...,+}(f_m(.,t))\in S^{(0,...,0,\frac{H}{2},\frac{H}{2})}_2H(\mathbb{R}^m).
\end{align*}
Now, we note that for any $\alpha=(\alpha_1,...,\alpha_m)$ such that $\alpha_i\geq 0$, $S^{(\alpha_1,...,\alpha_m)}_2H(\mathbb{R}^m)$ is a Hilbert space endowed with the inner product:
\begin{align*}
\forall f,g\in S^{(\alpha_1,...,\alpha_m)}_2H(\mathbb{R}^m)\ <f,g>=\int_{\mathbb{R}^m}\big(\prod_{i=1}^m(1+|\xi_i|^2)^{\alpha_i}\big)\mathcal{F}(f)(\xi)\overline{\mathcal{F}(g)(\xi)}d\xi
\end{align*}
Moreover, it is straightforward to check that $S^{(0,...,0,H/2,H/2)}_2H(\mathbb{R}^m)$ is isomorphically isometric to $L^2(\mathbb{R}^{m-2})\otimes S^{(H/2,H/2)}_2H(\mathbb{R}^2)$. The next lemma allows us to consider the restriction of function in $S^{(0,...,0,H/2,H/2)}_2H(\mathbb{R}^m)$ to the hyperplane $\Delta_m$.\\
\\
\textbf{Lemma 3.12}: There exists a linear and continuous operator $T$ from $S^{(0,...,0,H/2,H/2)}_2H(\mathbb{R}^m)$ to $S^{(0,...,0,H-1/2)}_2H(\mathbb{R}^{m-1})$ such that for any $f\in S^{(0,...,0,H/2,H/2)}_2H(\mathbb{R}^m)\cap C(\mathbb{R}^m)$, we have:
\begin{align*}
T(f)=f|_{\Delta_m}.
\end{align*}
\textbf{Proof}: We note $\tilde{T}$ the operator defined by:
\begin{align*}
\forall\psi\in L^2(\mathbb{R}^{m-2}), \forall\phi\in S_2^{(\frac{H}{2},\frac{H}{2})}H(\mathbb{R}^2),\ \tilde{T}(\psi\otimes\phi)=\psi\otimes Tr(\phi),
\end{align*}
and extended by linearity. Moreover, we have:
\begin{align*}
\|\tilde{T}(\psi\otimes\phi)\|_{L^2({\mathbb{R}^{m-2}})\otimes H^{H-\frac{1}{2}}(\mathbb{R})}\leq C_H\|\psi\otimes\phi\|_{L^2{\mathbb{R}^{m-2}}\otimes S_2^{(\frac{H}{2},\frac{H}{2})}H(\mathbb{R}^2)},
\end{align*} 
where $C_H$ is a strictly positive constant depending on $H$ uniquely. Thus, by a standard argument, we can extend $\tilde{T}$ to a continuous and linear operator from $L^2(\mathbb{R}^{m-2})\otimes S^{(H/2,H/2)}_2H(\mathbb{R}^2)$ to $L^2(\mathbb{R}^{m-2})\otimes H^{H-\frac{1}{2}}(\mathbb{R})$. We note this extension $\tilde{T}$ and it is exactly the tensor product operator $E\otimes Tr$ where $E$ is the identity operator from $L^2(\mathbb{R}^{m-2})$ to $L^2(\mathbb{R}^{m-2})$. Using the linear isometry between $L^2(\mathbb{R}^{m-2})\otimes S^{(H/2,H/2)}_2H(\mathbb{R}^2)$ and $S^{(0,...,0,H/2,H/2)}_2H(\mathbb{R}^m)$ (and between $L^2(\mathbb{R}^{m-2})\otimes H^{H-1/2}(\mathbb{R})$ and $S_2^{(0,...,0,H-1/2)}H(\mathbb{R}^{m-1})$ respectively), we obtain a linear continuous operator $T$ from $S^{(0,...,0,H/2,H/2)}_2H(\mathbb{R}^m)$ to $S_2^{(0,...,0,H-1/2)}H(\mathbb{R}^{m-1})$. It is straightforward to check that this operator reduces to the restriction operator over $\Delta_m$ for functions $f\in S^{(0,...,0,H/2,H/2)}_2H(\mathbb{R}^m)\cap C(\mathbb{R}^m)$ since $\Delta_m=\mathbb{R}^{m-2}\times\Delta_2$. $\Box$\\
\\
We end this preamble by discussing the hypothesis \textbf{H1}. First of all, a function in $S^{(0,...,0,H/2,H/2)}_2H(\mathbb{R}^m)$ admits an integral representation as a convolution of a function in $L^2(\mathbb{R}^m)$ with the Bessel polypotential kernel $G_{(0,...,0,H/2,H/2)}$ (see \cite{MR0192223} and chapter $24.12$ of \cite{samko1993fractional} for a definition of this kernel). Thus, for any $t\in (a,b)$ and $m\geq 2$, we have:
\begin{align*}
\exists g_{m}(.,t)\in L^2(\mathbb{R}^m),\ I^{0,...,0,\frac{H}{2},\frac{H}{2}}_{+,...,+}(f_m(,t))=G_{(0,...,0,H/2,H/2)}\ast g_m(,t).
\end{align*}
Using the Fourier operator, we have:
\begin{align*}
\mathcal{F}(f_m(,t))=\dfrac{(-i\xi_{m-1})^{\frac{H}{2}}(-i\xi_m)^{\frac{H}{2}}}{\prod_{i=m-1}^m (1+|\xi_i|^2)^{\frac{H}{4}}}\mathcal{F}(g_m(,t)).
\end{align*}
A direct application of Marcinkiewicz multiplier theorem (see \cite{MR0290095} chapter $4.3$ Theorem $3$ and chapter $4.6$) leads to the following inequality:
\begin{align*}
\|f_m(,t)\|_{L^2(\mathbb{R}^m)}\leq C_H\|g_m(,t)\|_{L^2(\mathbb{R}^m)}.
\end{align*}
Moreover, the $L^2(\mathbb{R}^m)$-norm of the function $g_m(,t)$ is equivalent to the norm of $I^{0,...,0,H/2,H/2}_{+,...,+}(f_m(,t))$ as an element of $S^{(0,...,0,H/2,H/2)}_2H(\mathbb{R}^m)$ (see \cite{MR0192223}). We are now ready to state the theorem regarding the variance of the stochastic integral of a random process with respect to the Rosenblatt process.\\
\\
\textbf{Theorem 3.13}: Let $(a,b)\subset\mathbb{R}_+$ with $0\leq a<b<\infty$. Let $\{\phi_t;t\in (a,b)\}$ be a stochastic process such that for all $t\in(a,b)$, $\phi_t\in(L^2)$. Moreover, let us assume that:
\begin{align*}
\forall t\in (a,b),\ \forall m\geq 2,\ f_m(.,t)\ verifies\ \mathbf{H1},\\ 
\sum_{m=0}^{+\infty}(m+2)!\bigg(\int_{(a,b)}\|g_m(.;t)\|^{q}_{L^2(\mathbb{R}^m)}dt\bigg)^{\frac{2}{q}}<+\infty,
\end{align*}
where $\phi_t=\sum_{m=0}^{+\infty}I_m(f_m(.,t))$ and $q=\frac{1}{H}(1+\epsilon)$ for some $\epsilon\in(\frac{1-H}{3H-1}\vee (2H-1);1)$.
Thus, we have:
\begin{align*}
&\mathbb{E}[(\int_{(a,b)}\phi_t\diamond\dot{X}^H_tdt)^2]=H(2H-1)\int_{(a,b)}\int_{(a,b)}|t-s|^{2(H-1)}\mathbb{E}[\phi_t\phi_s]dsdt\\&+4\sqrt{\frac{H(2H-1)}{2}}\int_{(a,b)}\int_{(a,b)}|t-s|^{H-1}\mathbb{E}[D_{\sqrt{d(H)}\delta_s\circ I^{\frac{H}{2}}_+}(\phi_t)D_{\sqrt{d(H)}
_t\circ I^{\frac{H}{2}}_+}(\phi_s)]dsdt\\&+\int_{(a,b)}\int_{(a,b)}\mathbb{E}[(D_{\sqrt{d(H)}\delta_s\circ I^{\frac{H}{2}}_+})^2(\phi_t)(D_{\sqrt{d(H)}\delta_t\circ I^{\frac{H}{2}}_+})^2(\phi_s)]dsdt.
\end{align*}
\textbf{Proof}: For any $t\in (a,b)$, by the Wiener-Itô Theorem, we have:
\begin{align*}
\phi_t=\sum_{m=0}^{+\infty}I_m(f_m(.,t))
\end{align*}
Let us fix $m$. By the definition of the operator $D^*_{\sqrt{d(H)}\delta_t\circ I^{\frac{H}{2}}_+}$ and by Pettis integration (in $(S)^*$ as well as in $S'(\mathbb{R})$), we obtain:
\begin{align*}
\int_{(a;b)}I_m(f_m(.,t))\diamond\dot{X}^H_tdt=d(H)I_{m+2}\left(\int_{(a,b)}(\delta_t\circ I^{\frac{H}{2}}_+)^{\otimes 2}\hat{\otimes} f_m(.,t)dt\right)
\end{align*}
Let $\xi\in S(\mathbb{R})$. We have:
\begin{align*}
<&\int_{(a,b)}(\delta_t\circ I^{\frac{H}{2}}_+)^{\otimes 2}\hat{\otimes} f_m(.,t)dt;\xi^{\otimes(m+2)}>=\int_{(a,b)}<(\delta_t\circ I^{\frac{H}{2}}_+)^{\otimes 2}\hat{\otimes} f_m(.,t);\xi^{\otimes(m+2)}>dt\\
=&\int_{\mathbb{R}^{m+2}}\bigg(\int_{(a,b)}\dfrac{(t-x_1)^{\frac{H}{2}-1}_+}{\Gamma(\frac{H}{2})}\dfrac{(t-x_2)^{\frac{H}{2}-1}_+}{\Gamma(\frac{H}{2})}f_m(x_3,...,x_{m+2},t)dt\bigg)\xi^{\otimes(m+2)}(x_1,...,x_{m+2})dx_1...dx_{m+2}
\end{align*}
Moreover, using Lemma 2.3, we have:
\begin{align*}
(d(H))^2||\int_{(a,b)}\dfrac{(t-.)^{\frac{H}{2}-1}_+}{\Gamma(\frac{H}{2})}\dfrac{(t-.)^{\frac{H}{2}-1}_+}{\Gamma(\frac{H}{2})}f_m(.,t)dt||^2_{L^2(\mathbb{R}^{m+2})}=\frac{H(2H-1)}{2}\\\int_{(a,b)}\int_{(a,b)}|t-s|^{2(H-1)}<f_m(.,t);f_m(.,s)>_{L^2(\mathbb{R}^m)}dtds
\end{align*}
which is finite by Hardy-Littlewood-Sobolev inequality (see Theorem 4.3 in \cite{MR1817225}) with $n=1$, $\lambda=2-2H$, $p=r=\frac{1}{H}$. Indeed, 
\begin{align*}
|\int_{(a,b)}\int_{(a,b)}|t-s|^{2(H-1)}<f_m(.,t);f_m(.,s)>_{L^2(\mathbb{R}^m)}dtds|\leq
&\int_{(a,b)}\int_{(a,b)}|t-s|^{2(H-1)}\|f_m(.,t)\|_{L^2(\mathbb{R}^m)}\\
&\times\|f_m(.,s)\|_{L^2(\mathbb{R}^m)}dtds,\\
&\leq C_H\bigg(\int_{(a,b)}\|f_m(.,t)\|^{\frac{1}{H}}_{L^2(\mathbb{R}^m)}dt\bigg)^{2H},\\
&\leq C_{H,a,b,q}\bigg(\int_{(a,b)}\|f_m(.,t)\|^{q}_{L^2(\mathbb{R}^m)}dt\bigg)^{\frac{2}{q}},\\
&\leq C_{H,a,b,q}\bigg(\int_{(a,b)}\|g_m(.,t)\|^{q}_{L^2(\mathbb{R}^m)}dt\bigg)^{\frac{2}{q}}<\infty.
\end{align*}
Thus,
\begin{align*}
d(H)\int_{(a,b)}(\delta_t\circ I^{\frac{H}{2}}_+)^{\otimes 2}\hat{\otimes} f_m(.,t)dt=d(H)sym\left(\int_{(a,b)}\dfrac{(t-x_1)^{\frac{H}{2}-1}_+}{\Gamma(\frac{H}{2})}\dfrac{(t-x_2)^{\frac{H}{2}-1}_+}{\Gamma(\frac{H}{2})}f_m(x_3,...,x_{m+2},t)dt\right)
\end{align*}
where $sym()$ denotes the symmetrization's operation and we have the following estimate ensuring that $\int_{(a,b)}\phi_t\diamond\dot{X}^H_tdt$ is in $(L)^2$:
\begin{align*}
\mathbb{E}[|\int_{(a,b)}\phi_t\diamond\dot{X}^H_tdt|^2]&\leq \frac{H(2H-1)}{2} \sum_{m=0}^{+\infty}(m+2)!\int_{(a,b)}\int_{(a,b)}|t-s|^{2(H-1)}<f_m(.,t);f_m(.,s)>_{L^2(\mathbb{R}^m)}dtds,\\
&\leq \frac{H(2H-1)}{2} \sum_{m=0}^{+\infty}(m+2)!\int_{(a,b)}\int_{(a,b)}|t-s|^{2(H-1)}\|f_m(.,t)\|_{L^2(\mathbb{R}^m)}\|f_m(.,s)\|_{L^2(\mathbb{R}^m)}dtds,\\
&\leq C_H\frac{H(2H-1)}{2} \sum_{m=0}^{+\infty}(m+2)!\bigg(\int_{(a,b)}\|f_m(.,t)\|^{\frac{1}{H}}_{L^2(\mathbb{R}^m)}dt\bigg)^{2H},\\
&\leq C_{H,a,b,q}\frac{H(2H-1)}{2}\sum_{m=0}^{+\infty}(m+2)!\bigg(\int_{(a,b)}\|f_m(.,t)\|^{q}_{L^2(\mathbb{R}^m)}dt\bigg)^{\frac{2}{q}}<\infty.
\end{align*}
Now, we compute precisely the symmetrization of the previous function in order to obtain the exact formula for the variance of $\int_{(a,b)}\phi_t\diamond\dot{X}^H_tdt$. By definition, we have:
\begin{align*}
sym&\left(\int_{(a,b)}\dfrac{(t-x_1)^{\frac{H}{2}-1}_+}{\Gamma(\frac{H}{2})}\dfrac{(t-x_2)^{\frac{H}{2}-1}_+}{\Gamma(\frac{H}{2})}f_m(x_3,...,x_{m+2},t)dt\right)=\\&\frac{1}{(m+2)!}\int_{(a,b)}\sum_{\sigma\in \mathcal{S}_{m+2}}\left(\dfrac{(t-x_{\sigma(1)})^{\frac{H}{2}-1}_+}{\Gamma(\frac{H}{2})}\dfrac{(t-x_{\sigma(2)})^{\frac{H}{2}-1}_+}{\Gamma(\frac{H}{2})}f_m(x_{\sigma(3)},...,x_{\sigma(m+2)},t)\right)dt\\
&=\frac{1}{(m+2)!}\int_{(a,b)}\left(\sum_{j=1}^{m+2}\dfrac{(t-x_j)^{\frac{H}{2}-1}_+}{\Gamma(\frac{H}{2})}\sum_{\sigma\in\mathcal{S}^j_{m+1}}\dfrac{(t-x_{\sigma(1)})^{\frac{H}{2}-1}_+}{\Gamma(\frac{H}{2})}f_m(x_{\sigma(2)},...,x_{\sigma(m+1)},t)\right)dt
\end{align*}
where $\mathcal{S}^j_{m+1}$ is the set of permutations of the set $\{1,...,j-1,j+1,...,m+2\}$.
\begin{align*}
=\frac{1}{(m+2)!}\int_{(a,b)}\left(\sum_{j=1}^{m+2}\sum_{k=1,k\ne j}^{m+2}\dfrac{(t-x_j)^{\frac{H}{2}-1}_+}{\Gamma(\frac{H}{2})}\dfrac{(t-x_k)^{\frac{H}{2}-1}_+}{\Gamma(\frac{H}{2})}\sum_{\sigma\in\mathcal{S}^{j,k}_{m}}f_m(x_{\sigma(1)},...,x_{\sigma(m+2)},t)\right)dt
\end{align*}
where $\mathcal{S}^{j,k}_{m}$ is the set of permutations of the set $\{1,...,m+2\}\setminus\{j,k\}$. But, the kernel $f_m(.,t)$ is symmetric in the $m$ first coordinates. Thus,
\begin{align*}
=\frac{1}{(m+2)(m+1)}\int_{(a,b)}\left(\sum_{j=1}^{m+2}\sum_{k=1,k\ne j}^{m+2}\dfrac{(t-x_j)^{\frac{H}{2}-1}_+}{\Gamma(\frac{H}{2})}\dfrac{(t-x_k)^{\frac{H}{2}-1}_+}{\Gamma(\frac{H}{2})}f_m(\underline{\mathbf{x}}_{j,k},t)\right)dt
\end{align*}
where $\underline{\mathbf{x}}_{j,k}$ is the vector $(x_1,...,x_{m+2})$ without the coordinates $x_j$ and $x_k$. We can now compute the $L^2$-norm of the symmetrization. We obtain:
\begin{align*}
=\frac{1}{[(m+2)(m+1)]^2}\sum_{j=1}^{m+2}\sum_{k=1,k\ne j}^{m+2}\sum_{i=1}^{m+2}\sum_{l=1,l\ne i}^{m+2}\int_{(a,b)}\int_{(a,b)}\bigg(&\int_{\mathbb{R}^{m+2}}\dfrac{(t-x_j)^{\frac{H}{2}-1}_+}{\Gamma(\frac{H}{2})}\dfrac{(t-x_k)^{\frac{H}{2}-1}_+}{\Gamma(\frac{H}{2})}\\&\dfrac{(s-x_i)^{\frac{H}{2}-1}_+}{\Gamma(\frac{H}{2})}\dfrac{(s-x_l)^{\frac{H}{2}-1}_+}{\Gamma(\frac{H}{2})}f_m(\underline{\mathbf{x}}_{j,k},t)
f_m(\underline{\mathbf{x}}_{i,l},s)d\underline{\mathbf{x}}\bigg)dtds
\end{align*}
To pursue, we need to enumerate the different values of the integral over $\mathbb{R}^{m+2}$ and to count the number of times such values appear. Since $j\ne k$ and $i\ne l$, we have at most equality between two couples of subscripts, namely, $i=j$ and $l=k$ or $j=l$ and $i=k$ which leads to $2(m+2)(m+1)$ possibilities. On the contrary, we can have every subscripts distinct corresponding to $(m+2)(m+1)m(m-1)$ possibilities. Finally, the $4(m+2)(m+1)m$ last possibilities correspond to the case where two subscripts are equal and all the others disinct. These three cases lead to three different values for the integral over $\mathbb{R}^{m+2}$. Specifically, we obtain:
\begin{align*}
\int_{\mathbb{R}^{m+2}}\dfrac{(t-x_j)^{\frac{H}{2}-1}_+}{\Gamma(\frac{H}{2})}\dfrac{(t-x_k)^{\frac{H}{2}-1}_+}{\Gamma(\frac{H}{2})}\dfrac{(s-x_i)^{\frac{H}{2}-1}_+}{\Gamma(\frac{H}{2})}\dfrac{(s-x_l)^{\frac{H}{2}-1}_+}{\Gamma(\frac{H}{2})}f_m(\underline{\mathbf{x}}_{j,k},t)
f_m(\underline{\mathbf{x}}_{i,l},s)d\underline{\mathbf{x}}=\\
\begin{cases}
\dfrac{(\beta(1-H,\frac{H}{2}))^2}{(\Gamma(\frac{H}{2}))^4}|t-s|^{2(H-1)}<f_m(,t);f_m(,s)>_{L^2(\mathbb{R}^{m})}\\
\dfrac{\beta(1-H,\frac{H}{2})}{(\Gamma(\frac{H}{2}))^2}|t-s|^{H-1}<(\delta_s\circ I^{\frac{H}{2}}_+)\otimes_1f_m(,t);(\delta_t\circ I^{\frac{H}{2}}_+)\otimes_1f_m(,s)>_{L^2(\mathbb{R}^{m-1})}\\
<(\delta_s\circ I^{\frac{H}{2}}_+)\otimes_1((\delta_s\circ I^{\frac{H}{2}}_+)\otimes_1f_m(,t));(\delta_t\circ I^{\frac{H}{2}}_+)\otimes_1((\delta_t\circ I^{\frac{H}{2}}_+)\otimes_1f_m(,s))>_{L^2(\mathbb{R}^{m-2})}
\end{cases}
\end{align*}
And, we obtain for the $L^2$ norm of the symmetrization:
\begin{align*}
&\frac{1}{(m+1)(m+2)}\bigg[2\dfrac{(\beta(1-H,\frac{H}{2}))^2}{(\Gamma(\frac{H}{2}))^4}\int_{(a,b)}\int_{(a,b)}|t-s|^{2(H-1)}<f_m(,t);f_m(,s)>_{L^2(\mathbb{R}^{m})}dsdt\\&+4m\dfrac{\beta(1-H,\frac{H}{2})}{(\Gamma(\frac{H}{2}))^2}\int_{(a,b)}\int_{(a,b)}|t-s|^{H-1}<(\delta_s\circ I^{\frac{H}{2}}_+)\otimes_1f_m(,t);(\delta_t\circ I^{\frac{H}{2}}_+)\otimes_1f_m(,s)>_{L^2(\mathbb{R}^{m-1})}dsdt\\&+m(m-1)\int_{(a,b)}\int_{(a,b)}<(\delta_s\circ I^{\frac{H}{2}}_+)\otimes_1((\delta_s\circ I^{\frac{H}{2}}_+)\otimes_1f_m(,t));(\delta_t\circ I^{\frac{H}{2}}_+)\otimes_1((\delta_t\circ I^{\frac{H}{2}}_+)\otimes_1f_m(,s))>_{L^2(\mathbb{R}^{m-2})}dsdt\bigg]
\end{align*}
We want to bound the three former terms in order to invert the sum and the integrals and to obtain the final formula for the variance. We denote the three integrals respectively by $(I)$, $(II)$ and $(III)$. The first one, $(I)$, has already been dealt with. Let us bound the second one:
\begin{align*}
(II)=\int_{(a,b)}\int_{(a,b)}|t-s|^{H-1}<I_{+,...,+}^{(0,..,0,\frac{H}{2})}(f_m(,t))(,s);I_{+,...,+}^{(0,..,0,\frac{H}{2})}(f_m(,s))(,t)>_{L^2(\mathbb{R}^{m-1})}dsdt,
\end{align*}
By Cauchy-Schwarz inequality, we obtain:
\begin{align*}
|(II)|\leq \int_{(a,b)}\int_{(a,b)}|t-s|^{H-1}\|I_{+,...,+}^{(0,..,0,H/2)}(f_m(,t))(,s)\|_{L^2(\mathbb{R}^{m-1})}\|I_{+,...,+}^{(0,..,0,H/2)}(f_m(,s))(,t)\|_{L^2(\mathbb{R}^{m-1})}dsdt.
\end{align*}
First of all, by Theorem $24.1$ of \cite{samko1993fractional} and since $\bigg(\int_{(a,b)}\|f_m(.;t)\|^{q}_{L^2(\mathbb{R}^m)}dt\bigg)^{1/q}<\infty$, we note that:
\begin{align*}
\bigg(\int_{(a,b)}\bigg(\int_{\mathbb{R}}\|I_{+,...,+}^{(0,..,0,\frac{H}{2})}(f_m(,t))(,s)\|_{L^2(\mathbb{R}^{m-1})}^{\frac{2}{1-H}}ds\bigg)^{\frac{q(1-H)}{2}}dt\bigg)^{\frac{1}{q}}<\infty.
\end{align*}
Or stated otherwise, the measurable function $(s,t)\rightarrow \|I_{+,...,+}^{(0,..,0,H/2)}(f_m(,t))(,s)\|_{L^2(\mathbb{R}^{m-1})}\mathbf{1}_{(a,b)}(t)$ is an element of $L^{(2/(1-H),q)}(\mathbb{R}^2)$. Similarly, using Theorem $2.51$ of \cite{MR2424078}, the measurable function $(s,t)\rightarrow \|I_{+,...,+}^{(0,..,0,H/2)}(f_m(,s))(,t)\|_{L^2(\mathbb{R}^{m-1})}\mathbf{1}_{(a,b)}(s)$ is an element of $L^{(q,2/(1-H))}(\mathbb{R}^2)$. Then, using a generalized version of Hölder inequality in the framework of mixed-norm spaces (see Theorem $2.49$ of \cite{MR2424078}), the measurable function $(s,t)\rightarrow \|I_{+,...,+}^{(0,..,0,H/2)}(f_m(,t))(,s)\|_{L^2(\mathbb{R}^{m-1})}\times\|I_{+,...,+}^{(0,..,0,H/2)}(f_m(,s))(,t)\|_{L^2(\mathbb{R}^{m-1})}\times\mathbf{1}_{(a,b)\times(a,b)}(s,t)$ is an element of $L^{(r,r)}(\mathbb{R}^2)$ where $1/r=(1-H)/2+1/q$. Moreover, for $\epsilon\in((1-H)/(3H-1)\vee (2H-1);1)$, we have $r>1/H$. Thus, by Hölder inequality with $p^*=r$, we obtain:
\begin{align*}
|(II)|\leq\bigg(\int_{(a,b)\times (a,b)}|t-s|^{p(H-1)}dsdt\bigg)^{\frac{1}{p}}\bigg(\int_{(a,b)\times (a,b)}\|I_{+,...,+}^{(0,..,0,\frac{H}{2})}(f_m(,t))(,s)\|^{p^*}_{L^2(\mathbb{R}^{m-1})}\\
\times\|I_{+,...,+}^{(0,..,0,\frac{H}{2})}(f_m(,s))(,t)\|^{p^*}_{L^2(\mathbb{R}^{m-1})}dsdt\bigg)^{\frac{1}{p^*}}.
\end{align*}
where,
\begin{align*}
\forall p\in (1,\frac{1}{1-H})\quad\int_{(a,b)\times (a,b)}|t-s|^{p(H-1)}dsdt<\infty.
\end{align*}
Using the generalized version of Hölder inequality for mixed-norm spaces, the permutation inequality for mixed norms, and Theorem $24.1$ of \cite{samko1993fractional}, we get:
\begin{align*}
|(II)|\leq C_{a,b,H,\epsilon} \bigg(\int_{(a,b)}\|f_m(;t)\|^{q}_{L^2(\mathbb{R}^{m})}dt\bigg)^{\frac{2}{q}}.
\end{align*}
In order to conclude, we need to bound the last integral, $(III)$. We have
\begin{align*}
(III)= \int_{(a,b)}\int_{(a,b)}<T\big(I^{(0,...,0,\frac{H}{2},\frac{H}{2})}_{+,...,+}(f_m(,t))\big)(,s);T\big(I^{(0,...,0,\frac{H}{2},\frac{H}{2})}_{+,...,+}(f_m(,s))\big)(,t)>_{L^2(\mathbb{R}^{m-2})}dsdt.
\end{align*}
By Cauchy-Schwarz inequality, we obtain:
\begin{align*}
|(III)|\leq \int_{(a,b)}\int_{(a,b)}\|T\big(I^{(0,...,0,\frac{H}{2},\frac{H}{2})}_{+,...,+}(f_m(,t))\big)(,s)\|_{L^2(\mathbb{R}^{m-2})}\times\|T\big(I^{(0,...,0,\frac{H}{2},\frac{H}{2})}_{+,...,+}(f_m(,s))\big)(,t)\|_{L^2(\mathbb{R}^{m-2})}dsdt.
\end{align*}
By Lemma 3.12, we have:
\begin{align*}
\|T\big(I^{(0,...,0,\frac{H}{2},\frac{H}{2})}_{+,...,+}(f_m(,t))\big)\|_{L^2(\mathbb{R}^{m-1})}&\leq \|T\big(I^{(0,...,0,\frac{H}{2},\frac{H}{2})}_{+,...,+}(f_m(,t))\big)\|_{S_2^{(0,...,0,H-1/2)}H(\mathbb{R}^{m-1})}\\
&\leq C_H \|I^{(0,...,0,\frac{H}{2},\frac{H}{2})}_{+,...,+}(f_m(,t))\|_{S_2^{(0,...,0,H/2,H/2)}H(\mathbb{R}^m)}\\
&\leq C_H'\|g_m(,t)\|_{L^2(\mathbb{R}^m)}.
\end{align*}
Thus,
\begin{align*}
\int_{(a,b)}\bigg(\int_{\mathbb{R}}\|T\big(I^{(0,...,0,\frac{H}{2},\frac{H}{2})}_{+,...,+}(f_m(,t))\big)(,s)\|^2_{L^2(\mathbb{R}^{m-2})}ds\bigg)^{\frac{q}{2}}dt<\infty.
\end{align*}
Since $q=1/H(1+\epsilon)$ with $\epsilon> (1-H)/(3H-1)\vee (2H-1)$, we have $q>2$. Thus,
\begin{align*}
\big(\int_{(a,b)\times (a,b)}\|T\big(I^{(0,...,0,\frac{H}{2},\frac{H}{2})}_{+,...,+}(f_m(,t))\big)(,s)\|^2_{L^2(\mathbb{R}^{m-2})}dsdt\big)^{\frac{1}{2}}\leq C_{a,b,H,\epsilon}\big(\int_{(a,b)}\|T\big(I^{(0,...,0,\frac{H}{2},\frac{H}{2})}_{+,...,+}(f_m(,t))\big)\|^q_{L^2(\mathbb{R}^{m-1})}dt\big)^{\frac{1}{q}},
\end{align*}
for some $C_{a,b,H,\epsilon}>0$. Using Cauchy-Schwarz inequality, we obtain:
\begin{align*}
|(III)|\leq C_{a,b,H,\epsilon} \bigg(\int_{(a,b)}\|g_m(;t)\|^{q}_{L^2(\mathbb{R}^{m})}dt\bigg)^{\frac{2}{q}}.
\end{align*}
Multiplying by $(m+2)!(d(H))^2$ and using the condition:
\begin{align*}
\sum_{m=0}^{+\infty}(m+2)!\bigg(\int_{(a,b)}\|g_m(.;t)\|^{q}_{L^2(\mathbb{R}^m)}dt\bigg)^{\frac{2}{q}}<+\infty,
\end{align*}
we can invert the sum and the integrals. Using the definition of $D_{\sqrt{d(H)}\delta_t\circ I^{\frac{H}{2}}_+}$ as well as the isometry for multiple Wiener-Itô integrals, we get the following formula for the $(L^2)$ norm of $\int_{(a,b)}\phi_t\diamond\dot{X}^H_tdt$.
\begin{align*}
&\mathbb{E}[(\int_{(a,b)}\phi_t\diamond\dot{X}^H_tdt)^2]=H(2H-1)\int_{(a,b)}\int_{(a,b)}|t-s|^{2(H-1)}\mathbb{E}[\phi_t\phi_s]dsdt\\&+4\sqrt{\frac{H(2H-1)}{2}}\int_{(a,b)}\int_{(a,b)}|t-s|^{H-1}\mathbb{E}[D_{\sqrt{d(H)}\delta_s\circ I^{\frac{H}{2}}_+}(\phi_t)D_{\sqrt{d(H)}\delta_t\circ I^{\frac{H}{2}}_+}(\phi_s)]dsdt\\&+\int_{(a,b)}\int_{(a,b)}\mathbb{E}[(D_{\sqrt{d(H)}\delta_s\circ I^{\frac{H}{2}}_+})^2(\phi_t)(D_{\sqrt{d(H)}\delta_t\circ I^{\frac{H}{2}}_+})^2(\phi_s)]dsdt.\quad\Box
\end{align*}
\\
\textbf{Remark 3.14}:
\begin{itemize}
\item These conditions are quite similar to the ones of Theorem $13.16$ in the case of the white noise integral (see chapter 13 of \cite{MR1387829}).
\item The first term appearing in the formula for the variance is classical and is linked to the space of distribution for which we can define the Wiener integrals with respect to the Rosenblatt process (see section 3 of \cite{MR2374640}). Actually, since the Rosenblatt process and fractional Brownian motion possess the same covariance structure, this space of distribution is the same than the one for fBm (see \cite{MR1790083} for further details).
\item We note the appearance of a second trace term which differs significantly from the Gaussian case.
\end{itemize}

In the rest of this section, we derive an Itô formula for a certain class of functionals of the Rosenblatt process. In the framework of white noise distribution theory, generalizations of Itô formula have been obtained for Brownian motion (See chapter 13.5 in \cite{MR1387829}), for fractional Brownian motion (see \cite{MR1956473}), and for convoluted Lévy processes (see \cite{MR2544099}). When one wants to derive an Itô formula with respect to a stochastic process, the first step is to identify for which class of functionals this formula will hold. At least, we want the polynomials to belong to this class. So, in the sequel, using the multiplicative formula from Malliavin calculus, we get an Itô formula for $x^2$ and $x^3$ in the $(S)^*$ sense. It appears that the series expansion intervening in the analytic representation around the origin of the logarithm of the translated characteristic function of the Rosenblatt process at time $t$ plays a fundamental role. Then, using Schwartz's Paley-Wiener theorem, we identify the class of functionals for which this Itô formula is true.\\
\\
\textbf{Lemma 3.15}: Let $(p,q)\in\mathbb{N}^*$. Let $f\in \hat{L}^2(\mathbb{R}^p)$ and $g\in\hat{L}^2(\mathbb{R}^q)$. Then, we have:
$$I_p(f)I_q(g)=\sum_{r=0}^{p\wedge q}r!C_r^pC_r^qI_{p+q-2r}(f\otimes_r g)$$
where $C_r^p=\frac{p!}{r!(p-r)!}$ and $\otimes_r$ denotes the contraction of order $r$.\\
\\
\textbf{Proof}: See proposition $1.1.3$ of \cite{MR2200233}.$\Box$\\
\\
\textbf{Theorem 3.16}: Let $(a,b)>0$ such that $a<b$. Then in $(S)^*$:
$$(X^H_b)^2-(X^H_a)^2=2\int_a^bX^H_sdX^H_s+b^{2H}-a^{2H}+4\int_a^bdX^{H,2}_s$$
and
$$(X^H_b)^3-(X^H_a)^3=3\int_a^b(X^H_s)^2dX^H_s+6H\int_a^bs^{2H-1}X^H_sds+12\int_a^bX^H_sdX^{H,2}_s+\kappa_3(X^H_1)(b^{3H}-a^{3H})+24\int_a^bdX^{H,3}_s.$$
\textbf{Proof}: Let $t>0$ and $\xi\in S(\mathbb{R})$. Using Lemma 3.13, we get:
$$(X^H_t)^2=I_4(f^H_t\hat{\otimes}f^H_t)+4I_2(f^H_t\otimes_1 f^H_t)+2||f^H_t||_{L^2(\mathbb{R}^2)}^2.$$
and,
$$(X^H_t)^3=I_6(f^H_t\hat{\otimes}f^H_t\hat{\otimes}f^H_t)+8I_4((f^H_t\hat{\otimes}f^H_t)\otimes_1 f^H_t)+4I_4((f^H_t\otimes_1f^H_t)\otimes f^H_t)+12 I_2((f^H_t\hat{\otimes} f^H_t)\otimes_2f^H_t)$$
$$+16 I_2((f^H_t\otimes_1f^H_t)\otimes_1 f^H_t)+2||f^H_t||^2_{L^2(\mathbb{R}^2)}I_2(f^H_t)+8<f^H_t\otimes_1 f^H_t;f^H_t>.$$
Then, using the definition of the $S$-transform, we have:
$$S((X^H_t)^2)(\xi)=(<f^H_t;\xi^{\otimes 2}>)^2+4<f^H_t\otimes_1 f^H_t;\xi^{\otimes 2}>+2||f^H_t||_{L^2(\mathbb{R}^2)}^2$$
$$S((X^H_t)^3)(\xi)=(<f^H_t;\xi^{\otimes 2}>)^3+8<(f^H_t\hat{\otimes}f^H_t)\otimes_1 f^H_t,\xi^{\otimes 4}>+4<(f^H_t\otimes_1f^H_t)\otimes f^H_t;\xi^{\otimes 4}>+$$
$$12<(f^H_t\hat{\otimes} f^H_t)\otimes_2f^H_t;\xi^{\otimes 2}>+16<(f^H_t\otimes_1f^H_t)\otimes_1 f^H_t;\xi^{\otimes 2}>+2||f^H_t||^2_{L^2(\mathbb{R}^2)}<f^H_t;\xi^{\otimes 2}>+8<f^H_t\otimes_1 f^H_t;f^H_t>.$$
By (2) and the cumulant's formula, we have:
$$8<f^H_t\otimes_1 f^H_t;f^H_t>=8\sum_{n=1}^{+\infty}(\lambda_nt^H)^3=8t^{3H}\dfrac{\kappa_3(X^H_1)}{2^22!}=t^{3H}\kappa_3(X^H_1).$$
and,
$$2||f^H_t||^2_{L^2(\mathbb{R})}=t^{2H}$$
Moreover, we have:
$$4<(f^H_t\otimes_1f^H_t)\otimes f^H_t;\xi^{\otimes 4}>=4<f^H_t\otimes_1f^H_t;\xi^{\otimes 2}><f^H_t;\xi^{\otimes 2}>$$
$$8<(f^H_t\hat{\otimes}f^H_t)\otimes_1 f^H_t,\xi^{\otimes 4}>=8<f^H_t;\xi^{\otimes 2}>||T_t(\xi)||^2_{L^2(\mathbb{R})}$$
Finally, we get, by Fubini's theorem:
$$12<(f^H_t\hat{\otimes} f^H_t)\otimes_2f^H_t;\xi^{\otimes 2}>=12<<f^H_t\hat{\otimes} f^H_t;\xi^{\otimes 2}>;f^H_t>.$$
By definition of the symmetric tensor product:
$$<f^H_t\hat{\otimes} f^H_t;\xi^{\otimes 2}>=\dfrac{1}{4!}\sum_{\sigma\in \mathcal{S}_4}\int_{\mathbb{R}^2}f^H_t(x_{\sigma(1)},x_{\sigma(2)})f^H_t(x_{\sigma(3)},x_{\sigma(4)})\xi^{\otimes 2}(x_1,x_2)dx_1dx_2$$
The integral inside the sum can only have two values (depending on $\sigma$):
$$\int_{\mathbb{R}^2}f^H_t(x_{\sigma(1)},x_{\sigma(2)})f^H_t(x_{\sigma(3)},x_{\sigma(4)})\xi^{\otimes 2}(x_1,x_2)dx_1dx_2=
\begin{cases}
<f^H_t;\xi^{\otimes 2}>f^H_t & \sigma=(\sigma',\sigma'') \\
T_t(\xi)^{\otimes 2} & \sigma\ne (\sigma',\sigma'')
\end{cases}$$
where $(\sigma',\sigma'')\in (S_2(\{1,2\}^{\{1,2\}})\times S_2(\{3,4\}^{\{3,4\}}))\cup (S_2(\{1,2\}^{\{3,4\}})\times S_2(\{3,4\}^{\{1,2\}}))$. Thus, we get:
$$<f^H_t\hat{\otimes} f^H_t;\xi^{\otimes 2}>=\dfrac{1}{3}<f^H_t;\xi^{\otimes 2}>f^H_t+\dfrac{2}{3}T_t(\xi)^{\otimes 2}$$
To conclude, we have:
$$12<(f^H_t\hat{\otimes} f^H_t)\otimes_2f^H_t;\xi^{\otimes 2}>=4<f^H_t;\xi^{\otimes 2}>||f^H_t||^2_{L^2(\mathbb{R}^2)}+8<f^H_t;T_t(\xi)^{\otimes 2}>$$
Once again, using (2), we obtain:
$$<f^H_t;T_t(\xi)^{\otimes 2}>=<(f^H_t\otimes_1 f^H_t)\otimes_1 f^H_t;\xi^{\otimes 2}>$$
And so, we obtain the following expression for the S-transform of $(X^H_t)^3$:
$$S((X^H_t)^3)(\xi)=(<f^H_t;\xi^{\otimes 2}>)^3+12<f^H_t;\xi^{\otimes 2}>||T_t(\xi)||^2_{L^2(\mathbb{R})}+3t^{2H}<f^H_t;\xi^{\otimes 2}>$$
$$+24<(f^H_t\otimes_1 f^H_t)\otimes_1 f^H_t;\xi^{\otimes 2}>+t^{3H}\kappa_3(X^H_1)$$
By differentiation, we get:
$$\dfrac{d}{dt}\left[S((X^H_t)^3)(\xi)\right]=(3(<f^H_t;\xi^{\otimes 2}>)^2+12||T_t(\xi)||^2_{L^2(\mathbb{R})}+3t^{2H})S(\dot{X}_t^H)(\xi)$$
$$+(12\dfrac{d}{dt}[||T_t(\xi)||^2_{L^2(\mathbb{R})}]+6Ht^{2H-1})S(X_t^H)(\xi)+24\dfrac{d}{dt}[<(f^H_t\otimes_1 f^H_t)\otimes_1 f^H_t;\xi^{\otimes 2}>]$$
$$+3Ht^{3H-1}\kappa_3(X^H_1)$$
Using the expression of the S-transform of $(X^H_t)^2$:
$$\dfrac{d}{dt}\left[S((X^H_t)^3)(\xi)\right]=3S((X^H_t)^2\diamond \dot{X}_t^H)(\xi)+6(Ht^{2H-1}+2\dfrac{d}{dt}[||T_t(\xi)||^2_{L^2(\mathbb{R})}])S(X_t^H)(\xi)$$ 
$$+6(H\dfrac{t^{3H-1}}{2}\kappa_3(X^H_1)+2^2\dfrac{d}{dt}[<(f^H_t\otimes_1 f^H_t)\otimes_1 f^H_t;\xi^{\otimes 2}>])$$
And by Lemma 3.7:
$$\dfrac{d}{dt}\left[S((X^H_t)^3)(\xi)\right]=3S((X^H_t)^2\diamond \dot{X}_t^H)(\xi)+(6Ht^{2H-1}+12S(\dot{X}^{H,2}_t)(\xi))S(X_t^H)(\xi)$$ 
$$+3Ht^{3H-1}\kappa_3(X^H_1)+24S(\dot{X}^{H,3}_t)(\xi)$$
In order to conclude, we have to show the $(S)^*$-integrability of the distribution processes $(X^H_t)^2\diamond \dot{X}_t^H$, $t^{2H-1}X_t^H$ and $X_t^H\diamond\dot{X}^{H,2}_t$ on $(a,b)$. By definition of the Wick product, we have:
$$S((X^H_t)^2\diamond\dot{X}^H_t)(\xi)=S((X^H_t)^2)(\xi)S(\dot{X}^H_t)(\xi)$$
Using the Chaos expansion of $(X^H_t)^2$ and Lemma 3.4,
$$S((X^H_t)^2\diamond\dot{X}^H_t)(\xi)=d(H)(I^{\frac{H}{2}}_+(\xi)(t))^2\left((S(X^H_t)(\xi))^2+4<f^H_t\otimes_1 f^H_t;\xi^{\otimes 2}>+t^{2H}\right)$$
And by Lemma 2.4,
\begin{align*}
S((X^H_t)^2\diamond\dot{X}^H_t)(\xi)=&d(H)(I^{\frac{H}{2}}_+(\xi)(t))^2\bigg[\big(d(H)\int_{0}^t[I_+^{\frac{H}{2}}(\xi)(s)]^2ds\big)^2\\
&+4d(H)\sqrt{\frac{H(2H-1)}{2}}
\int_{0}^t\int_{0}^tI_+^{\frac{H}{2}}(\xi)(s)I_+^{\frac{H}{2}}(\xi)(r)|s-r|^{H-1}dsdr+t^{2H}\bigg]
\end{align*}
Since $I^{\frac{H}{2}}_+(\xi)\in C^{\infty}(\mathbb{R})$, $t\rightarrow S((X^H_t)^2\diamond\dot{X}^H_t)(\xi)$ is continuous on $(a,b)$ and then measurable.\\
Moreover using Lemma 3.3 and the fact that $\int_0^1\int_0^1|s-r|^{H-1}dsdr<\infty$, we get easily the last point of theorem 3.9 ensuring the $(S)^*$-integrability of the distribution process $(X^H_t)^2\diamond\dot{X}^H_t$. We proceed as previously to get the $(S)^*$-integrability of the distribution processes $X_t^H\diamond\dot{X}^{H,2}_t$ and $t^{2H-1}X_t^H$.$\Box$\\
\\
\textbf{Remark 3.17}:
\begin{itemize} 
\item For $x^3$, we note the appearance of two terms involving the third derivative of this monomial. As it is written in \cite{MR2374640}, one cannot hope to obtain an Itô formula which ends with a second derivative term for the Rosenblatt process.
\item Before stating the main results, we need an estimate of the radius of convergence of the power series expansion of the logarithm of the translated characteristic function of the Rosenblatt process at time $t$ in a neighborhood of the origin. Using,
$|c_k|\leq (\frac{1}{H(2H-1)})^{\frac{k}{2}}$ and Theorem 2.9, the power series expansion converge if $|\sqrt{2}\theta t^H|<1$. For $t\in (a,b)$ where $0<a\leq b<\infty$ we obtain the sufficient condition $|\theta|<\frac{1}{\sqrt{2}b^H}$.
\end{itemize}

\textbf{Theorem 3.18}: Let $(a,b)\in \mathbb{R}^*_+$ such that $a\leq b<\infty$. Let $F$ be an entire analytic function of the complex variable verifying:
$$\forall N\in\mathbb{N},\exists C_N>0,\forall z\in\mathbb{C}\quad |F(z)|\leq C_N\dfrac{\exp(\frac{1}{\sqrt{2}b^H}|\Im(z)|)}{(1+|z|)^N}$$
Then, we have in $(S)^*$:
$$F(X^H_b)-F(X^H_a)=\int_a^bF^{(1)}(X^H_t)\diamond\dot{X}^H_tdt+\sum_{k=2}^{\infty}\left(H\kappa_{k}(X^H_1)\int_a^b\dfrac{t^{Hk-1}}{(k-1)!}F^{(k)}(X_t^H)dt+2^k\int_a^bF^{(k)}(X_t^H)\diamond\dot{X}^{H,k}_tdt\right)$$
\textbf{Proof}: Let $(a,b)\in \mathbb{R}^*_+$ such that $a\leq b<\infty$ and let $F$ be as in the theorem. First of all, let's compute the $S$-transform of $F(X^H_t)$ for $t\in (a,b)$. By definition, we have:
$$\forall\xi\in S(\mathbb{R})\quad S(F(X^H_t))(\xi)=\mathbb{E}[F(X^H_t):\exp(<;\xi>):]=\int_{S'(\mathbb{R})}F((X^H_t)(x))d\mu_{\xi}(x)$$
$$=\int_{S'(\mathbb{R})}F((X^H_t)(x+\xi))d\mu(x)$$
Using the Transfer theorem and the Parseval-Plancherel's theorem, we obtain:
$$S(F(X^H_t))(\xi)=\dfrac{1}{2\pi}\int_{\mathbb{R}}\mathcal{F}(F)(\theta)\mathbb{E}^{\mu_{\xi}}[\exp(i\theta X^H_t)]d\theta$$
Moreover, by the Paley-Wiener theorem (see \cite{MR0493420}, Theorem $9.11$), the Fourier transform of $F$ is in $C_0^{\infty}(\mathbb{R})$ and its support is contained in the ball $\{\theta:|\theta|\leq\frac{1}{\sqrt{2}b^H}\}$. Thus,
\begin{align*}
S(F(X^H_t))(\xi)=\dfrac{1}{2\pi}\int_{supp(\mathcal{F}(F))}\mathcal{F}(F)(\theta)&\exp(i\theta<f^H_t;\xi^{\otimes 2}>)\exp\bigg(\frac{1}{2}\sum_{k=2}^{+\infty}\frac{(2i\theta t^H)^k}{k}\Big(\sqrt{\frac{(2H-1)H}{2}}\Big)^kc_k\bigg)\\
&\times\exp\bigg(\sum_{k=2}^{+\infty}(2i\theta)^k||T^{\frac{k}{2}}_t(\xi)||_2^2\bigg)d\theta
\end{align*}
where, we have (by Theorem 2.9):
$$||T^{\frac{k}{2}}_t(\xi)||_2^2=<\underbrace{(...((f_t^H\otimes_1f_t^H)\otimes_1f_t^H)...\otimes_1f_t^H)}_{k-1 \times\otimes_1};\xi^{\otimes 2}>=S(X^{H,k}_t)(\xi)$$
We differentiate with respect to $t$ and get:
\begin{align*}
\dfrac{d}{dt}S(F(X^H_t))(\xi)=\dfrac{1}{2\pi}\int_{supp(\mathcal{F}(F))}\mathcal{F}(F)(\theta)&\dfrac{d}{dt}\bigg[\exp(i\theta<f^H_t;\xi^{\otimes 2}>)\exp\bigg(\frac{1}{2}\sum_{k=2}^{+\infty}\frac{(2i\theta t^H)^k}{k}\Big(\sqrt{\frac{(2H-1)H}{2}}\Big)^kc_k\bigg)\\
&\times\exp\bigg(\sum_{k=2}^{+\infty}(2i\theta)^kS(X^{H,k}_t)(\xi)\bigg)\bigg]d\theta.
\end{align*}
The differentiation of the translated characteristic function with respect to $t\in (a,b)$ for $\theta\in supp(\mathcal{F}(F))$ gives:
$$\dfrac{d}{dt}\left[\mathbb{E}^{\mu_{\xi}}[\exp(i\theta X^H_t)]\right]=(I)+(II)+(III)$$
where,
$$(I)=i\theta S(\dot{X}^H_t)(\xi)\mathbb{E}^{\mu_{\xi}}[\exp(i\theta X^H_t)],$$
$$(II)=\dfrac{H}{2}\sum_{k=2}^{+\infty}\bigg[t^{Hk-1}(2i\theta)^k\Big(\sqrt{\frac{(2H-1)H}{2}}\Big)^kc_k\bigg]\mathbb{E}^{\mu_{\xi}}[\exp(i\theta X^H_t)],$$
and by Lemma 3.7,
$$(III)=\sum_{k=2}^{+\infty}\bigg[(2i\theta)^kS(\dot{X}^{H,k}_t)(\xi)\bigg]\mathbb{E}^{\mu_{\xi}}[\exp(i\theta X^H_t)].$$
Since we integrate over $supp(\mathcal{F}(F))$ strictly contained in $\{\theta:|\theta|\leq\frac{1}{\sqrt{2}b^H}\}$ where the power series are absolutely convergent, we have:
\begin{align*}
\dfrac{d}{dt}S(F(X^H_t))(\xi)=&S(\dot{X}^H_t)(\xi)\bigg(\dfrac{1}{2\pi}\int_{supp(\mathcal{F}(F))}\mathcal{F}(F)(\theta)i\theta\mathbb{E}^{\mu_{\xi}}[\exp(i\theta X^H_t)]d\theta\bigg)\\
+&\sum_{k=2}^{\infty}\bigg[2^{k-1}Ht^{Hk-1}\Big(\sqrt{\frac{(2H-1)H}{2}}\Big)^kc_k\Big(\dfrac{1}{2\pi}\int_{supp(\mathcal{F}(F))}\mathcal{F}(F)(\theta)(i\theta)^k\mathbb{E}^{\mu_{\xi}}[\exp(i\theta X^H_t)]d\theta\Big)\bigg]\\
+&\sum_{k=2}^{\infty}\bigg[2^kS(\dot{X}^{H,k}_t)(\xi)\Big(\dfrac{1}{2\pi}\int_{supp(\mathcal{F}(F))}\mathcal{F}(F)(\theta)(i\theta)^k\mathbb{E}^{\mu_{\xi}}[\exp(i\theta X^H_t)]d\theta\Big)\bigg]
\end{align*}
Using $(i\theta)^k\mathcal{F}(F)(\theta)=\mathcal{F}(F^{(k)})(\theta)$, the formula for $S(F(X^H_t))(\xi)$ and $(3)$, we obtain:
\begin{align*}
\dfrac{d}{dt}S(F(X^H_t))(\xi)=S(\dot{X}^H_t)(\xi)S(F^{(1)}(X^H_t))(\xi)&+\sum_{k=2}^{\infty}\bigg[H\dfrac{t^{Hk-1}}{(k-1)!}\kappa_k(X^H_1)S(F^{(k)}(X^H_t))(\xi)\bigg]\\
+&\sum_{k=2}^{\infty}\bigg[2^kS(\dot{X}^{H,k}_t)(\xi)S(F^{(k)}(X^H_t))(\xi)\bigg]
\end{align*}
To conclude we need to show the $(S)^*$-integrability over $(a,b)$ of the three former terms using Theorem 3.9. From the formula of Lemma 3.7 we get the following estimate:
\begin{align}
|S(\dot{X}^{H,k}_t)(\xi)|\leq C_1(H)\sqrt{\frac{H(2H-1)}{2}}^{k}\Big[&(H(k-1)+1)t^{H(k-1)}(||I_+^{\frac{H}{2}}(\xi)||_{\infty})^2\int_0^1\int_0^1K^{k-2}_1(u,v)dudv\\\nonumber
&+2t^{H(k-1)+1}||I_+^{\frac{H}{2}}(\xi)||_{\infty}||I_+^{\frac{H}{2}}(\xi^{(1)})||_{\infty}\int_0^1\int_0^1K^{k-2}_1(u,v)dudv\Big].
\end{align}
By Cauchy-Schwartz inequality, we have:
$$\int_0^1\int_0^1K^{k-2}_1(u,v)dudv\leq \left(\dfrac{1}{2}\dfrac{2}{(2H-1)H}\right)^{\frac{k}{2}}<\infty$$
Moreover, we have:
$$|S(F^{(k)}(X^H_t))(\xi)|\leq \dfrac{1}{2\pi}\int_{supp(\mathcal{F}(F))}|\mathcal{F}(F)(\theta)||\theta|^kd\theta$$
Writting $supp(\mathcal{F}(F))=[m_{F};M_{F}]\subset \{\theta:|\theta|<\frac{1}{\sqrt{2}b^H}\}$:
$$|S(F^{(k)}(X^H_t))(\xi)|\leq \dfrac{1}{2\pi}||\mathcal{F}(F)||_{L^1(\mathbb{R})}\max\{|m_F|^k;|M_{F}|^k\}$$
So the general term of the third series is bounded by:
\begin{align*}
|2^kS(\dot{X}^{H,k}_t)(\xi)S(F^{(k)}(X^H_t))(\xi)|\leq& C_2(H)\dfrac{1}{2\pi}(\sqrt{2})^{k}||\mathcal{F}(F)||_{L^1(\mathbb{R})}\max\{|m_F|^k;|M_{F}|^k\}\\
&\Big[(H(k-1)+1)t^{H(k-1)}(||\xi||_{\infty}+||\xi^{(1)}||_{\infty}+||\xi||_1)^2\\
+&2t^{H(k-1)+1}(||\xi||_{\infty}+||\xi^{(1)}||_{\infty}+||\xi||_1)(||\xi^{(1)}||_{\infty}+||\xi^{(2)}||_{\infty}+||\xi^{(1)}||_1)\Big].
\end{align*}
Thus, this term is $(S)^*$-integrable over $(a,b)$ and we have:
\begin{align*}
\int_a^b|2^kS(\dot{X}^{H,k}_t)(\xi)&S(F^{(k)}(X^H_t))(\xi)|dt\leq C_2(H)\dfrac{1}{2\pi}(\sqrt{2})^{k}||\mathcal{F}(F)||_{L^1(\mathbb{R})}\max\{|m_F|^k;|M_{F}|^k\}\\
&\Big[(b^{H(k-1)+1}-a^{H(k-1)+1})(||\xi||_{\infty}+||\xi^{(1)}||_{\infty}+||\xi||_1)^2\\
+&\frac{2}{H(k-1)+2}(b^{H(k-1)+2}-a^{H(k-1)+2})(||\xi||_{\infty}+||\xi^{(1)}||_{\infty}+||\xi||_1)(||\xi^{(1)}||_{\infty}+||\xi^{(2)}||_{\infty}+||\xi^{(1)}||_1)\Big].
\end{align*}
Since $\max\{|m_F|;|M_{F}|\}<\frac{1}{\sqrt{2}b^H}$, we can invert the sum and the integral. We proceed as previously to prove the $(S)^*$-integrability of the two other terms and we obtain the result.$\Box$\\
\\
\textbf{Remark 3.19}:
\begin{itemize}
\item The Itô formula obtained stresses the fundamental role of the analytic expansion around the origin of the logarithm of the translated characteristic function. In particular, the non-zero cumulants of the Rosenblatt distribution are crucial.   
\item Moreover, this result confirms the Itô formula obtained for the monomial $x^2$ and $x^3$. Nevertheless, the polynomials do not belong to the class of functionals for which the Itô formula holds. Thus, we need a generalization of the Paley-Wiener theorem in order to get the result for a bigger class which contains the polynomials. This is the aim of the next theorem.
\end{itemize}
\textbf{Theorem 3.20}: Let $(a,b)\in \mathbb{R}^*_+$ such that $a\leq b<\infty$. Let $F$ be an entire analytic function of the complex variable verifying:
$$\exists N\in\mathbb{N},\exists C>0,\forall z\in\mathbb{C}\quad |F(z)|\leq C(1+|z|)^N\exp(\frac{1}{\sqrt{2}b^H}|\Im(z)|)$$
Then, we have in $(S)^*$:
$$F(X^H_b)-F(X^H_a)=\int_a^bF^{(1)}(X^H_t)\diamond\dot{X}^H_tdt+\sum_{k=2}^{\infty}\left(H\kappa_{k}(X^H_1)\int_a^b\dfrac{t^{Hk-1}}{(k-1)!}F^{(k)}(X_t^H)dt+2^k\int_a^bF^{(k)}(X_t^H)\diamond\dot{X}^{H,k}_tdt\right)$$
\textbf{Proof}: Let $(a,b)\in \mathbb{R}^*_+$ such that $a\leq b<\infty$ and let $F$ be as in the theorem. As previously, we have to compute the $S$-transform of $F(X^H_t)$ for any $t\in(a,b)$. Using Cauchy-Schwartz inequality and the definition of $F$, we have:
$$\forall\xi\in S(\mathbb{R})\quad |S(F(X^H_t))(\xi)|\leq \mathbb{E}^{\mu_{\xi}}[|F(X^H_t)|]\leq \mathbb{E}[|F(X^H_t)|^2]^{\frac{1}{2}}\mathbb{E}[(:\exp(<;\xi>):)^2]^{\frac{1}{2}}$$
$$|S(F(X^H_t))(\xi)|\leq C\mathbb{E}[(1+|X^H_t|)^{2N}]^{\frac{1}{2}}\exp(\frac{||\xi||_2^2}{2})<\infty$$
since $X^H_t$ is in the second Wiener's chaos and, by hypercontractivity (see Corollary 2.8.14 in \cite{MR2962301}), has finite moments of any order. Thus, we get:
$$S(F(X^H_t))(\xi)=\dfrac{1}{2\pi}<\mathcal{F}(F);\mathbb{E}^{\mu_{\xi}}[\exp(i. X^H_t)]>$$
where $<;>$ has to be understood as the duality bracket since, by Schwartz's Paley-Wiener theorem (see \cite{MR0493420}, Theorem $9.12$), $\mathcal{F}(F)$ is in $S'(\mathbb{R})$, has a compact support contained in the ball $\{\theta:|\theta|\leq \frac{1}{\sqrt{2}b^H}\}$ and is of order $N$. Let $[m_F;M_F]$ be such that: 
$$supp(\mathcal{F}(F))\subsetneq [m_F;M_F]\subsetneq \{\theta:|\theta|\leq\frac{1}{\sqrt{2}b^H}\}.$$
Since $\mathcal{F}(F)$ is a distribution with compact support and since $\mathbb{E}^{\mu_{\xi}}[\exp(i. X^H_t)]\in C^{\infty}(\mathbb{R})$ (by Theorem 2.9), we have:
$$S(F(X^H_t))(\xi)=\dfrac{1}{2\pi}<\mathcal{F}(F);\psi_F(.)\mathbb{E}^{\mu_{\xi}}[\exp(i. X^H_t)]>$$
where $\psi_F\in C^{\infty}_0(\mathbb{R})$, is equal to $1$ in a neighborhood of $supp(\mathcal{F})$ and $supp(\psi_F)=[m_F;M_F]$. We differentiate with respect to $t$ and get as previously:
$$\dfrac{d}{dt}S(F(X^H_t))(\xi)=\dfrac{1}{2\pi}<\mathcal{F}(F);\psi_F(.)\dfrac{d}{dt}[\mathbb{E}^{\mu_{\xi}}[\exp(i. X^H_t)]]>$$
with,
$$\dfrac{d}{dt}[\mathbb{E}^{\mu_{\xi}}[\exp(i\theta X^H_t)]]=(I)+(II)+(III)$$
where,
$$(I)=i\theta S(\dot{X}^H_t)(\xi)\mathbb{E}^{\mu_{\xi}}[\exp(i\theta X^H_t)],$$
$$(II)=\dfrac{H}{2}\sum_{k=2}^{+\infty}\bigg[t^{Hk-1}(2i\theta)^k\Big(\sqrt{\frac{(2H-1)H}{2}}\Big)^kc_k\bigg]\mathbb{E}^{\mu_{\xi}}[\exp(i\theta X^H_t)],$$
$$(III)=\sum_{k=2}^{+\infty}\bigg[(2i\theta)^kS(\dot{X}^{H,k}_t)(\xi)\bigg]\mathbb{E}^{\mu_{\xi}}[\exp(i\theta X^H_t)].$$
In order to invert the sum and the duality bracket, we need to estimate $<\mathcal{F}(F);\psi_F(.)(i.)^k\mathbb{E}^{\mu_{\xi}}[\exp(i. X^H_t)]>$ for large values of $k$ ($\geq N$). Since $N$ is the order of $\mathcal{F}(F)$, we have (see chapter 2.3 of \cite{MR1996773}):
$$|<\mathcal{F}(F);\psi_F(.)(i.)^k\mathbb{E}^{\mu_{\xi}}[\exp(i. X^H_t)]>|\leq C_F\sum_{l=0}^N\underset{\theta\in[m_F;M_F]}{\sup}|\frac{d^l}{d\theta^l}\left((i\theta)^k\mathbb{E}^{\mu_{\xi}}[\exp(i\theta X^H_t)]\right)|$$ 
where $C_F>0$. Using Leibniz's formula, we obtain:
\begin{align*}|<\mathcal{F}(F);\psi_F(.)(i.)^k&\mathbb{E}^{\mu_{\xi}}[\exp(i. X^H_t)]>|\leq C_F\sum_{l=0}^N\sum_{m=0}^lC^l_m\\&\times\max\{|m_F|^{k-m};|M_F|^{k-m}\}k(k-1)...(k-m+1)\mathbb{E}^{\mu_{\xi}}[|X^H_t|^{l-m}]
\end{align*}
Thus, using $c_k\leq (\frac{1}{H(2H-1)})^{\frac{k}{2}}$, we obtain the following estimate:
\begin{align*}
t^{Hk-1}2^k\Big(\sqrt{\frac{(2H-1)H}{2}}\Big)^k&c_k|<\mathcal{F}(F);\psi_F(.)(i.)^k\mathbb{E}^{\mu_{\xi}}[\exp(i. X^H_t)]>|\leq C_F\sum_{l=0}^N\sum_{m=0}^lC^l_m\\&\times(\sqrt{2})^kt^{Hk-1}\max\{|m_F|^{k-m};|M_F|^{k-m}\}k(k-1)...(k-m+1)\mathbb{E}^{\mu_{\xi}}[|X^H_t|^{l-m}]
\end{align*}
Since $\max\{|m_F|;|M_F|\}<\frac{1}{\sqrt{2}b^H}$, we can invert the sum and the duality bracket for $(II)$. Similarly, using (9), we can invert the sum and the duality bracket for $(III)$. We get:
\begin{align*}
\dfrac{d}{dt}S(F(X^H_t))(\xi)=S(\dot{X}^H_t)(\xi)S(F^{(1)}(X^H_t))(\xi)&+\sum_{k=2}^{\infty}\bigg[H\dfrac{t^{Hk-1}}{(k-1)!}\kappa_k(X^H_1)S(F^{(k)}(X^H_t))(\xi)\bigg]\\
+&\sum_{k=2}^{\infty}\bigg[2^kS(\dot{X}^{H,k}_t)(\xi)S(F^{(k)}(X^H_t))(\xi)\bigg].
\end{align*}
To prove the $(S)^*$-integrability of the three former terms, we proceed as for the proof of theorem 3.13 using the appropriate estimates.$\Box$\\
\\
\textbf{Remark 3.21}:
\begin{itemize}
\item Every polynomial verifies the condition of the theorem. Moreover, for such functionals, the sums are finite. 
\item Using the growth property of $F$ and the hypercontractivity of the second Wiener chaos elements, we can show that the left-hand side of the equality is in $(L^2)$.
\end{itemize}

\section{Comparison with other approaches}
In \cite{MR2374640}, the author introduces two definitions of the stochastic integral with respect to the Rosenblatt process. The fisrt one is a pathwise integral using the method of \cite{russo1993forward} and the second one is a Skorohod integral using elements of the Malliavin calculus. The Skorohod approach relies on the finite interval representation of the Rosenblatt process. From this representation by means of the double Skorohod integral with respect to Brownian motion, the author defines the stochastic integral with respect to the Rosenblatt process. In this section, we show that this Skorohod integral coincides with the Rosenblatt noise integral defined using the finite interval representation.\\
\\
\textbf{Definition-Theorem 4.20}: The Rosenblatt process $\{X^H_t\}_{t\in\mathbb{R}}$ is equal in distribution to the process $\{Z^H_t\}_{t\in [0,T]}$ defined by:
\begin{align*}
\forall t\in[0,T]\quad Z^H_t=c(H)\int_{[0;t]^2}\int_0^t\prod_{j=1}^2(\dfrac{s}{x_j})^{\frac{H}{2}}(s-x_j)_{+}^{\frac{H}{2}-1}dsdB_{x_1}dB_{x_2}
\end{align*}
\textbf{Proof}: See theorem 1.1 in \cite{pipiras2010regularization} or proposition 1 in \cite{MR2374640}. $\Box$\\
\\
\textbf{Remark 4.21}: It is important to note that even if the two processes have the same finite-dimensional distributions under $\mu$, this fact is not true anymore under the probability $\mu_{\xi}$ for some $\xi\in S(\mathbb{R})$. Indeed, one can show that for $\xi\in C^{\infty}_0(\mathbb{R},\mathbb{R}_+)$ such that $supp(\xi)\subset (0,+\infty)$:
\begin{align*}
\forall t\in[0;T]\quad S(X^H_t)(\xi)\ne S(Z^H_t)(\xi)
\end{align*}
\\
Using this representation, we can compute the $S$-transform of the process $\{Z^H_t\}$ as well as differentiate it in $(S)^*$.\\
\\
\textbf{Proposition 4.22}: The process $\{Z^H_t\}$ is $(S)^*$-differentiable and the $S$-trasnform of its derivative is equal to:
\begin{align*}
\forall\xi\in S(\mathbb{R})\forall t\in[0,T]\quad S(\dot{Z}^H_t)(\xi)=c(H)\big(\int_0^t\xi(x)(\frac{t}{x})^{\frac{H}{2}}(t-x)^{\frac{H}{2}-1}_{+}dx\big)^2
\end{align*}
\textbf{Proof}: First of all, we need to compute the $S$-transform of the process $\{Z^H_t\}$. By definition, we have:
\begin{align*}
\forall\xi\in S(\mathbb{R})\forall t\in[0,T]\quad S(Z^H_t)(\xi)=c(H)\int_{[0,t]^2}\Big(\int_0^t\prod_{j=1}^2(\dfrac{s}{x_j})^{\frac{H}{2}}(s-x_j)_{+}^{\frac{H}{2}-1}ds\Big)\xi(x_1)\xi(x_2)dx_1dx_2
\end{align*}
By Fubini's theorem, we have:
\begin{align*}
S(Z^H_t)(\xi)&=c(H)\int_0^t\Big(\int_0^t(\dfrac{s}{x})^{\frac{H}{2}}(s-x)_{+}^{\frac{H}{2}-1}\xi(x)dx\Big)^2ds\\
&=c(H)\int_0^t\Big(\int_0^s(\dfrac{s}{x})^{\frac{H}{2}}(s-x)_{+}^{\frac{H}{2}-1}\xi(x)dx\Big)^2ds
\end{align*}
By making the change of variable $x=us$, we obtain:
\begin{align*}
S(Z^H_t)(\xi)=c(H)\int_0^ts^{H}\Big(\int_0^1(1-u)^{\frac{H}{2}-1}_{+}\xi(su)u^{-\frac{H}{2}}du\Big)^2ds
\end{align*}
Since $\xi\in S(\mathbb{R})$ and $H\in(\frac{1}{2},1)$, we get:
\begin{align*}
\forall t>0\quad\dfrac{d}{dt}(S(Z^H_t)(\xi))&=c(H)t^{H}\Big(\int_0^1(1-u)^{\frac{H}{2}-1}_{+}\xi(tu)u^{-\frac{H}{2}}du\Big)^2\\
&=c(H)\Big(t^{\frac{H}{2}}\int_0^tx^{-\frac{H}{2}}(t-x)_{+}^{\frac{H}{2}-1}\xi(x)dx\Big)^2
\end{align*}
Moreover, for any $h>0$ such that $t+h\leq T$:
\begin{align}
\Big|\dfrac{S(Z^H_{t+h})(\xi)-S(Z^H_t)(\xi)}{h}\Big|&=c(H)\frac{1}{h}\Big|\int_t^{t+h}s^{H}\Big(\int_0^1(1-u)^{\frac{H}{2}-1}_{+}\xi(su)u^{-\frac{H}{2}}du\Big)^2ds\Big|\nonumber\\
&\leq c(H)T^{H}\|\xi\|^2_{\infty}\Big(\int_0^1(1-u)^{\frac{H}{2}-1}_{+}u^{-\frac{H}{2}}du\Big)^2<\infty
\end{align}
Thus, using Theorem 1.3, we obtain the result. $\Box$\\
\\
Similarly to the previous section, we can define a stochastic integral with respect to $\{Z^H_t\}$, using the Wick product with the stochastic distribution process $\{\dot{Z}^H_t\}$ and Pettis integration over $[0,T]$.\\
\\
\textbf{Definition-Theorem 4.23}: Let $\{\phi_t;t\in [0,T]\}$ be a $(S)^*$ stochastic process satisfying the assumptions of Theorem 3.9. Then, $\phi_t\diamond \dot{Z}^H_t$ is $(S)^*$ integrable over $[0,T]$ and:
\begin{align*}
\int_{[0,T]}\phi_tdZ^H_t=\int_{[0,T]}\phi_t\diamond \dot{Z}^H_tdt
\end{align*}
\textbf{Proof}: Using estimate $(10)$ and Theorem 3.9, one can show that $\phi_t\diamond \dot{Z}^H_t$ is $(S)^*$-integrable over $[0,T]$. $\Box$\\
\\
In the next proposition, we prove that this integral is an extension of the Skorohod integral defined in \cite{MR2374640}. For this purpose, we need to introduce elements of the Malliavin calculus with respect to Brownian motion. Let $\mathcal{S}$ be the class of smooth random variables of the form:
\begin{align*}
\forall(t_1,...,t_n)\in[0,T]\quad F=f(B_{t_1},...,B_{t_n})
\end{align*}
where $f\in\mathcal{C}^{\infty}_b(\mathbb{R}^n)$. Then, we define the Malliavin derivative on the class $\mathcal{S}$ by:
\begin{align*}
\forall t\in[0,T]\quad D_t(F)=\sum_{j=1}^{n}\dfrac{\partial f}{\partial x_j}(B_{t_1},...,B_{t_n})\mathbf{1}_{[0,t_j]}(t)
\end{align*}
This operator is unbounded, closable and can be extended to the closure of $\mathcal{S}$, denoted $\mathbb{D}^{k,p}$, with respect to the norm (see section 1.2 in \cite{MR2200233}):
\begin{align*}
\forall F\in\mathcal{S},k\geq 1,p\geq 2\quad ||F||^p_{k,p}=\mathbb{E}[|F|^p]+\sum_{j=1}^k\mathbb{E}[||D^j(F)||_{L^2([0,T]^j)}]
\end{align*}
where $D^{j}$ is the $j$th iteration of the operator $D$. For any $j\geq 1$, $D^{j}$ takes values in $L^2(\Omega; L^2([0,T]^j))$. We denote by $\delta$ the adjoint operator of $D$. It is an unbounded operator from a domain in $L^2(\Omega;L^2([0,T]))$ to $L^2(\Omega)$. An element $\phi\in L^2(\Omega;L^2([0,T]))$ belongs to the domain of $\delta$ if and only if:
\begin{align*}
\forall F\in\mathbb{D}^{1,2}\quad |\mathbb{E}[<DF;\phi>_{L^2([0;T])}]|\leq C_\phi\sqrt{\mathbb{E}[F^2]}
\end{align*}
Moreover, we have the following duality relationship:
\begin{align*}
\forall F\in\mathbb{D}^{1,2}\quad\mathbb{E}[\delta(\phi)F]=\mathbb{E}[<DF;\phi>_{L^2([0;T])}]
\end{align*}
For any $j\geq 1$, this formula extends to the adjoint operator of $D^{j}$, denoted by $\delta^j$, and we have:
\begin{align*}
\forall \phi\in Dom \delta^j\quad\forall F\in\mathbb{D}^{1,2}\quad\mathbb{E}[\delta^j(\phi)F]=\mathbb{E}[<D^{j}F,\phi>_{L^2([0;T]^j)}]
\end{align*}
For further properties regarding the multiple integral $\delta^j$, we refer the reader to \cite{david1988generalized}. We are now ready to state the definition of the Skorohod type integral with to respect to the Rosenblatt process introduced in \cite{MR2374640}.\\
\\
\textbf{Definition 4.24}: Let $\{\phi_t:t\in[0;T]\}$ be a square integrable stochastic process. Then, we define its Skorohod integral with respect to $\{Z_t\}_{t\in[0;T]}$ by:
\begin{align*}
\int_{0}^T\phi_t\delta Z^H_t= \delta^{2}(I_H(\phi))
\end{align*}
where $I_H(\phi)$ is defined by:
\begin{align*}
I_H(\phi)(s_1,s_2)=c(H)\int_{s_1\vee s_2}^T\phi_u(\dfrac{u}{s_1})^{\frac{H}{2}}(u-s_1)^{\frac{H}{2}-1}(\dfrac{u}{s_2})^{\frac{H}{2}}(u-s_2)^{\frac{H}{2}-1}du
\end{align*}
We state the main result of this section which links the two definitions of integral.\\
\\
\textbf{Proposition 4.25}: Let $\{\phi_t;t\in[0;T]\}$ be a stochastic process such that $\phi\in L^2(\Omega;\mathcal{H})\cap L^2([0,T];\mathbb{D}^{2,2})$ and $\mathbb{E}[\int_0^T\int_0^T||D^{2}_{s_1,s_2}\phi||^2_{\mathcal{H}}ds_1ds_2]<\infty$ where
\begin{align*}
\mathcal{H}=\{f:[0;T]\rightarrow\mathbb{R};\int_0^T\int_0^Tf(s)f(t)|t-s|^{2H-2}dsdt<\infty\}.
\end{align*}
Then, $\{\phi_t\}$ is Skorohod integrable and $(S)^*$-integrable with respect to the Rosenblatt process, $\{Z^H_t\}_{t\in[0;T]}$, and we have:
\begin{align*}
\int_0^T\phi_t\delta Z^H_t=\int_0^T\phi_t\diamond\dot{Z}^H_tdt
\end{align*}
\textbf{Proof}: The fact that $\{\phi_t\}$ is Skorohod integrable with respect to $\{Z^H_t\}$ follows from Lemma 1, section 6 in \cite{MR2374640}. Let us show that $\phi_t\diamond\dot{Z}^H_t$ is $(S)^*$-integrable. By definition, we have:
\begin{align*}
\forall\xi\in S(\mathbb{R})\quad |S(\phi_t\diamond\dot{Z}^H_t)(\xi)|&=|S(\phi_t)(\xi)||S(\dot{Z}^H_t)(\xi)|\\
&\leq c(H)t^{H}\|\xi\|^2_{\infty}|\mathbb{E}^{\mu}[\phi_t:\exp(<;\xi>):]|\Big(\int_0^1(1-u)^{\frac{H}{2}-1}u^{-\frac{H}{2}}du\Big)^2
\end{align*}
Using Cauchy-Schwartz inequality, we obtain:
\begin{align*}
|\mathbb{E}^{\mu}[\phi_t:\exp(<;\xi>):]|&\leq (\mathbb{E}[(\phi_t)^2])^{\frac{1}{2}}(\mathbb{E}[(:\exp(<;\xi>):)^2])^{\frac{1}{2}}\\
&\leq (\mathbb{E}[(\phi_t)^2])^{\frac{1}{2}}\exp(\frac{\|\xi\|^2_2}{2})
\end{align*}
Thus, we get the following estimate:
\begin{align*}
|S(\phi_t\diamond\dot{Z}^H_t)(\xi)|\leq  c(H)t^{H}(\mathbb{E}[(\phi_t)^2])^{\frac{1}{2}}\|\xi\|^2_{\infty}\exp(\frac{\|\xi\|^2_2}{2})\Big(\int_0^1(1-u)^{\frac{H}{2}-1}u^{-\frac{H}{2}}du\Big)^2
\end{align*}
Since $\phi\in L^2([0,T];\mathbb{D}^{2,2})$, we have,
\begin{align*}
\int_0^Tt^{H}(\mathbb{E}[(\phi_t)^2])^{\frac{1}{2}}dt<+\infty
\end{align*}
So, we conclude by Theorem 3.9. We have to show that the two integrals are equal. Let us compute the $S$-transform of the Skorohod integral of $\{\phi_t\}$ with respect to $\{Z^H_t\}$. By definition,
\begin{align*}
\forall\xi\in S(\mathbb{R})\quad S(\int_0^T\phi_t\delta Z^H_t)(\xi)=\mathbb{E}^{\mu}[\delta^2(I_H(\phi)):\exp(<;\xi>):]
\end{align*}
By the duality formula, we obtain:
\begin{align*}
S(\int_0^T\phi_t\delta Z^H_t)(\xi)&=\mathbb{E}^{\mu}[<I_H(\phi);D^2(:\exp(<;\xi>):)>_{L^2([0;T]^2)}]\\
&=\mathbb{E}^{\mu}[<I_H(\phi);\xi^{\otimes 2}>_{L^2([0;T]^2)}:\exp(<;\xi>):]
\end{align*}
Moreover, for any $\psi$ such that $I_H(\psi)\in L^2([0;T]^2)$,
\begin{align*}
<I_H(\psi);\xi^{\otimes 2}>_{L^2([0;T])}&=c(H)\int_{[0;T]^2}\Big(\int_{s_1\vee s_2}^T\psi(u)\prod_{j=1}^2(\dfrac{u}{s_j})^{\frac{H}{2}}(u-s_j)^{\frac{H}{2}-1}du\Big)\xi^{\otimes 2}(s_1,s_2)ds_1ds_2\\
&=c(H)\int_0^T\psi(u)\Big(\int_0^u(\dfrac{u}{s})^{\frac{H}{2}}(u-s)^{\frac{H}{2}-1}\xi(s)ds\Big)^2du
\end{align*}
Thus, we get:
\begin{align*}
S(\int_0^T\phi_t\delta Z^H_t)(\xi)=\int_0^T\mathbb{E}^{\mu}[\phi_t:\exp(<;\xi>):]S(\dot{Z}^H_t)(\xi)dt=S(\int_0^T\phi_t\diamond\dot{Z}^H_tdt)(\xi)
\end{align*}
And so we are done. $\Box$\\

\def\cprime{$'$}

\end{document}